\documentclass[final,onefignum,onetabnum]{siamart190516}

\usepackage[utf8]{inputenc}
\usepackage[T1]{fontenc}
\usepackage{amsfonts}
\usepackage{amssymb}
\usepackage{mathtools}
\usepackage{booktabs}
\usepackage{tabularx}
\usepackage{nicefrac}
\usepackage{microtype}
\usepackage{comment}
\usepackage{algpseudocode}
\usepackage{todonotes}
\usepackage{capt-of}
\usepackage{float}
\graphicspath{{figures/}}
\definecolor{lightblue}{HTML}{a1b4c7}
\definecolor{orange}{HTML}{ea8810}
\definecolor{silver}{HTML}{b0aba8}
\definecolor{rust}{HTML}{b8420f}
\definecolor{seagreen}{HTML}{23553c}
\colorlet{lightsilver}{silver!30!white}
\colorlet{darkorange}{orange!85!black}
\colorlet{darksilver}{silver!85!black}
\colorlet{darklightblue}{lightblue!85!black}
\colorlet{darkrust}{rust!85!black}
\colorlet{darkseagreen}{seagreen!85!black}
\hypersetup{
  colorlinks=true,
  linkcolor=darkrust,
  citecolor=darklightblue,
  urlcolor=darksilver,
  pdftitle={Kernelized Stable Fluids with Sparse Inverse Cholesky Preconditioning},
  pdfauthor={Lennart Scholz, Aras Bacho, and Houman Owhadi}
}

\makeatletter
\renewcommand{\paragraph}{%
  \@startsection{paragraph}{4}%
  {\z@}{1.0ex \@plus .5ex \@minus .2ex}{-.7em}%
  {\normalfont\normalsize\bfseries}%
}
\makeatother

\headers{Kernelized Stable Fluids}{L. Scholz, A. Bacho, and H. Owhadi}
\title{Kernelized Stable Fluids for simulating physically accurate solutions}

\author{%
  Lennart\ Scholz\thanks{Leibniz University Hannover,
  \email{l.scholz@stud.uni-hannover.de}}
  \and
  Aras\ Bacho\thanks{Computing and Mathematical Sciences,
  California Institute of Technology, Pasadena, CA 91125,
  \email{bacho@caltech.edu}}
  \and
  Houman\ Owhadi\thanks{Computing and Mathematical Sciences,
  California Institute of Technology, Pasadena, CA 91125,
  \email{owhadi@caltech.edu}}
}

\begin{document}
\raggedbottom
\maketitle

\begin{abstract}
The Stable Fluids (SF) algorithm \cite{StableFluid} is a fractional step method for solving the incompressible Navier–Stokes (NS) equations in an unconditionally stable way. It was initially developed for visual effects and allows real-time simulations because of its low computational complexity. The method solves the advection step by using the method of characteristics which is implemented via linear interpolation. This results in excessive diffusion of vorticity and is one of the main reasons for the physical inaccuracy of the algorithm. We propose Kernelized Stable Fluids with Sparse Cholesky Preconditioning (KSF-Chol) in order to address this shortcoming and make the algorithm thereby suitable for engineering problems and scientific computing. The algorithm is modified by replacing the linear interpolation approach with kernel interpolation in the advection step. We accelerate the computationally expensive kernel interpolation by using sparse Cholesky preconditioning and a near-/far-field approach and optimize both methods for computations on a uniform Cartesian grid. KSF-Chol achieves per time step a near-linear computational complexity of $\mathcal{O}(N\log^d(N/\varepsilon_{\mathrm{chol}}))$ in space and time.
We test KSF-Chol against reference solvers and numerical reference data on multiple 2D and 3D benchmarks and demonstrate a significantly reduced numerical diffusion as well as a competitive accuracy concerning the computation of relevant flow quantities.
\end{abstract}

\section{Introduction}
The dynamics of viscous incompressible fluids are governed by the NS equations. Numerical approaches for solving these equations usually demand a discretization of the governing equations through Finite Volumes (FV), Finite Elements (FE) or Finite Difference (FD) methods \cite{ferziger2020computational}. To circumvent the ill-conditioning of these methods for large time steps and achieve unconditional stability Jos Stam introduced the Stable Fluids algorithm \cite{StableFluid}, which is a fractional step method and uses the method of characteristics for fluid advection.
During each time step, the Stable Fluids algorithm successively applies a
forcing, the solution operations for the advection and diffusion step and finally performs a projection in order to achieve a divergence free velocity field. The numerical stability of Stable Fluids comes at the cost
of numerical dissipation, which arises mainly during the advection step. The step is solved by tracing
characteristics backward in time and evaluating the velocity field at the
resulting off-grid points. These
evaluations are performed by linear interpolation on a regular Cartesian grid and are therefore a weak point of the solver's architecture that we want to improve. While dissipation of vorticity
can be often acceptable for visually plausible fluid animations, it limits the
applicability of the original Stable Fluids algorithm when quantitatively accurate simulations for scientific applications
are required.

The central idea of this work is to improve the accuracy of the
advection step while keeping the stable fractional-step structure of the Stable
Fluids algorithm. We introduce \emph{Kernelized Stable Fluids with Sparse Cholesky
Preconditioning} (KSF-Chol), in which the linear interpolation used for
the advection step is replaced by kernel interpolation. The objective is to reduce the numerical diffusion introduced
during the advection and thereby obtain a method that is more
suitable for quantitative and scientific flow simulations. Further modifications introduced in KSF-Chol lead to a reduced 
computational cost of the usually very expensive kernel interpolation.

For $N$ grid points, one kernel interpolation requires
the solution of a dense linear system
$\mathbf{K}\boldsymbol{\alpha}=\mathbf{u}$ and a subsequent evaluation of the
kernel expansion $f(\cdot)=\sum_{j=1}^{N} K(\cdot,X_j)\alpha_j$ at $N$ points in our case. A direct factorization of
$\mathbf{K}$ would require $\mathcal{O}(N^3)$ time and
$\mathcal{O}(N^2)$ memory and is therefore prohibitively expensive for large
CFD grids. KSF-Chol solves the linear systems with a
preconditioned conjugate gradient (PCG) method. We use a sparse
Cholesky preconditioner \cite{sparse_chol_leibler}. Its non-aggregated construction costs $\mathcal{O}(N\rho_{\mathrm{chol}}^{3d})$ time and $\mathcal{O}(N\rho_{\mathrm{chol}}^d)$ memory, and each application costs $\mathcal{O}(N\rho_{\mathrm{chol}}^d)$. The aggregated construction in \cite{sparse_chol_leibler} costs
$\mathcal{O}(N\rho_{\mathrm{chol}}^{2d})$ time and
$\mathcal{O}(N\rho_{\mathrm{chol}}^d)$ space,
but is not used here. For certain elliptic Green's-function kernels, \cite{sparse_chol_leibler} establishes that an $\varepsilon_{\mathrm{chol}}$-accurate kernel matrix approximation can be obtained with $\rho_{\mathrm{chol}}\approx\log(N/\varepsilon_{\mathrm{chol}})$. All our benchmarks instead use $\rho_{\mathrm{chol}}=7.5$. For a constant $\rho_{\mathrm{chol}}$ we observe higher PCG iteration counts on larger grids and attribute this to a reduced preconditioner accuracy. We observed a constant iteration number within each simulation. Since $\mathbf{K}$ remains unchanged in every solver step, the sparse Cholesky factorization is done once during an initialization step and the resulting sparse Cholesky factor $\mathbf{L}$ is reused in every solver step.

For a translation-invariant kernel (Gaussian, Mat\'ern),
additional efficiency optimization is possible when the input points of the kernel form a uniform Cartesian grid. In this case, the kernel matrix has a multilevel
Toeplitz structure and its matrix-vector products can be evaluated as discrete
convolutions using FFTs. This reduces the time complexity of the matrix-vector product in the PCG step from
$\mathcal{O}(N^2)$ to $\mathcal{O}(N\log N)$ and avoids storing the
dense matrix. This optimization is restricted to the equidistant, uniform Cartesian discretisation of the computational domain, it cannot in general be applied to arbitrary point sets or
meshes. Furthermore, the evaluation points of the kernel expansion themselves are
off-grid and do not satisfy this Toeplitz structure, their evaluation is computationally expensive if it is performed in the dense manner.

A dense kernel evaluation with $N$ targets requires
$\mathcal{O}(N^2)$ pairwise interactions. We again exploit the uniform Cartesian grid geometry in order to
reduce this cost and employ an exact near-field and compressed far-field approximation, motivated by fast
multipole methods and low-rank compression of kernel interactions
\cite{FongDarve2009,MessnerBramasCoulaudDarve2012}. Each off-grid target is decomposed
into a Cartesian base-cell index and a normalized cell-coordinate in
$[-1,1]^d$. The dependence on the continuous cell coordinate is represented
by Chebyshev polynomials, while the
compressed far-field representation consists of $R$ spatial modes with $R \ll N$. Since these modes
depend only on integer grid offsets, their action on the kernel
coefficients is again a Toeplitz convolution and can be evaluated by FFTs. The
near-field is treated exactly. The resulting off-grid kernel evaluation has
$\mathcal{O}(N\log N)$ time and $\mathcal{O}(N)$ space complexity, instead
of the $\mathcal{O}(N^2)$ arithmetic cost of a direct kernel interpolation.

With $\rho_{\mathrm{chol}}\approx\log(N/\varepsilon_{\mathrm{chol}})$ the cost of KSF-Chol per time step is $\mathcal{O}(N\log^d(N/\varepsilon_{\mathrm{chol}}))$ in time and space.

\subsection{Summary of contributions}

The contributions of this work are summarized as follows:

\begin{itemize}

\item \textbf{Kernelized semi-Lagrangian advection.}
We replace the linear interpolation in the semi-Lagrangian advection step by
\[
    \psi(U)(Y)
    =
    K(Y,X)K(X,X)^{-1}U,
    \qquad
    X,Y\in(\mathbb{R}^{d})^{N},
    \quad
    U\in\mathbb{R}^{N}.
\]

\item \textbf{Sparse preconditioned CG method.}
The kernel systems
\[
    K(X,X)\boldsymbol{\alpha}=U
\]
are solved by PCG with a sparse inverse Cholesky preconditioner
\cite{sparse_chol_leibler}. For translation-invariant kernels on a uniform
Cartesian grid, $K(X,X)$ is multilevel Toeplitz and its matrix--vector products can be
evaluated by FFTs in $\mathcal{O}(N\log N)$ time complexity with $\mathcal{O}(N)$ space complexity
\cite{ChanJin2007,RakhubaOseledets2015}.

\item \textbf{Near-/far-field off-grid evaluation.}
The dense off-grid evaluation
\[
    \psi(U)(Y)=K(Y,X)\boldsymbol{\alpha}
\]
is split into an exact near field and a compressed far field.
For fixed approximation parameters, the evaluation requires
$\mathcal{O}(N\log N)$ time and $\mathcal{O}(N)$ space instead of
$\mathcal{O}(N^2)$ arithmetic complexity
\cite{FongDarve2009,MessnerBramasCoulaudDarve2012}.

\item \textbf{Numerical validation.}
We compare KSF-Chol with the original Stable Fluids algorithm, pseudo-spectral and finite-difference
reference solvers, and numerical reference data and investigate vorticity conservation,
residuals, vortex-shedding statistics, energy
spectra, and pressure/separation diagnostics.

\end{itemize}

\subsection{Review of literature}

A major obstacle to large-scale kernel methods is that the associated kernel
matrices are generally dense.
\cite{sparse_chol_leibler} introduced an approximation for such matrices using sparse inverse Cholesky
factorization. The same principle has also been extended to kernel matrices
arising in the numerical solution of nonlinear PDEs \cite{msichol}.

Beyond interpolation accuracy, kernel methods offer a flexible framework
for incorporating physical constraints and problem specific structure.
Gaussian Process Hydrodynamics (GPH) uses physics-informed kernels to
approximate solutions of the Euler and Navier--Stokes equations, incorporating
incompressibility, boundary conditions, and multiscale
turbulence statistics \cite{owhadi2023gph}.
Empirical kernels constructed from PDE solution snapshots can further adapt
to anisotropy, sharp gradients, and boundary behavior \cite{Bacho2026KROM}. 
Their probabilistic Gaussian-process formulation also provides a basis for
uncertainty quantification \cite{OwhadiScovelSchaefer2019}
and data assimilation through conditioning on observations.
Related approaches have been applied to inverse PDE problems, including
the identification of unknown physical parameters
\cite{ChenHosseiniOwhadiStuart2021}.
The circulation-preserving fluid formulation of
\cite{PicTurb} was likewise introduced to address the numerical dissipation of
vorticity that can arise in classical Stable Fluids methods.

\medskip

\noindent

\subsection{Outline of the article}
The remainder of the article is organized as follows. In
Section~\ref{sec:ksf-chol}, we introduce the KSF-Chol framework, starting from
the original Stable Fluids method and subsequently presenting the kernel based
velocity field reconstruction, sparse Cholesky preconditioning, and the
implementation details of the algorithm. Section~\ref{sec:numerical-experiments} presents the
benchmark problems, the evaluation methodology, and the numerical
comparison of KSF-Chol with classical numerical solvers and numerical
reference data. Finally, Section~\ref{sec:conclusions} summarizes the main
findings and discusses possible directions for future work.
\section{Kernelized Stable Fluids with Sparse Cholesky Preconditioning}
\label{sec:ksf-chol}
This section introduces the numerical framework of KSF-Chol.
Section~\ref{sec:stable-fluids} first reviews the original Stable Fluids
algorithm.
Section~\ref{sec:kernel-methods-stable-fluids} then introduces the
kernel methods used for reconstructing the continuous velocity fields from discrete data.
Section~\ref{sec:sparse-cholesky-preconditioning} presents the sparse Cholesky
preconditioner that accelerates the PCG method. Section~\ref{sec:cheb-svd-fft-evaluation} introduces the near-/far-field approach that accelerates the off-grid kernel evaluation while exploiting the uniform Cartesian grid structure.
Finally, Section~\ref{sec:ksf-chol-implementation} combines these methods
into the complete KSF-Chol algorithm and describes its numerical
implementation.
\subsection{Stable Fluids}
\label{sec:stable-fluids}

The original SF algorithm \cite{StableFluid} as presented in \ref{alg:SF_algo} solves the incompressible Navier--Stokes equations of the form
\begin{align}
\nabla \cdot \mathbf{u} &= 0, \\
\frac{\partial \mathbf{u}}{\partial t}
&= -(\mathbf{u}\cdot\nabla)\mathbf{u}
- \nabla p
+ \nu \Delta \mathbf{u}
+ \mathbf{f},
\label{eq:inkomp_eq}
\end{align}
where $\mathbf{u}$ denotes the velocity field, with components
$\mathbf{u}=(u,v)$ in two dimensions and $\mathbf{u}=(u,v,w)$ in three dimensions.
The scalar $p$ is the pressure, $\nu$ is the kinematic viscosity, and
$\mathbf{f}$ represents an external body force.

\begin{algorithm}
\caption{Stable Fluids}
\label{alg:SF_algo}

\normalsize
\algrenewcommand\alglinenumber[1]{\normalsize #1}

\begin{algorithmic}[1]

\State \textbf{Input:} spatial discretization $N_x \times N_y$, uniform Cartesian grid $\mathbf{x}=\{\mathbf{x}_i\}_{i=1}^{M}$ with
$\mathbf{x}_i=(x_i,y_i)$ in 2D and $M=N_xN_y$ cell centers

\vspace{0.3em}
\Statex \textbf{Initialization:}

\State Initialize velocity $\mathbf{u}^0$

\vspace{0.3em}
\For{$n = 0,1,\dots,N_t-1$}

\Statex \textbf{Add forcing:}

\State $\mathbf{u}^{\mathrm{force}}
       = \mathbf{u}^n + \Delta t\,\mathbf{f}$

\vspace{0.3em}
\Statex \textbf{Advection:}

\State $\mathbf{x}_{\mathrm{mid}}
       = \mathbf{x}
       - \dfrac{\Delta t}{2}\,
       \mathbf{u}^{\mathrm{force}}(\mathbf{x})$
        \Comment{$\mathbf{u}^{\mathrm{force}}$ is known at the grid points $\mathbf{x}$}

\State $\mathbf{x}^{\ast}
       = \mathbf{x}
       - \Delta t\,
       \mathbf{u}^{\mathrm{force}}(\mathbf{x}_{\mathrm{mid}})$
       \Comment{interpolation of $\mathbf{u}^{\mathrm{force}}$ at off-grid points $\mathbf{x}_{\mathrm{mid}}$}

\State $\mathbf{u}^{\mathrm{adv}}(\mathbf{x})
       = \mathbf{u}^{\mathrm{force}}(\mathbf{x}^{\ast})$
       \Comment{interpolation of $\mathbf{u}^{\mathrm{force}}$ at off-grid points $\mathbf{x}^{\ast}$}

\vspace{0.3em}
\Statex \textbf{Diffusion:}

\State $(I - \nu \Delta t\,\Delta)\,
       \mathbf{u}^{\mathrm{diff}}
       = \mathbf{u}^{\mathrm{adv}}$
        \Comment{solved by linear \texttt{PyFishPack} solver}

\vspace{0.3em}
\Statex \textbf{Projection:}

\State $\Delta q
       = \nabla \cdot \mathbf{u}^{\mathrm{diff}}$
       \Comment{solved by linear \texttt{PyFishPack} solver}

\State $\mathbf{u}^{n+1}
       = \mathbf{u}^{\mathrm{diff}} - \nabla q$

\EndFor

\vspace{0.3em}
\State \textbf{Return:} velocity fields $\{\mathbf{u}^n(\mathbf{x})\}_{n=0}^{N_t}$

\end{algorithmic}
\end{algorithm}

\subsection{Kernel Methods for Stable Fluids}
\label{sec:kernel-methods-stable-fluids}

Let $d$ denote the spatial dimension of the computational domain
$\Omega\subset\mathbb{R}^{d}$. The semi-Lagrangian advection step requires
evaluations of the velocity field (see interpolation of $\mathbf{u}^{\mathrm{force}}$ in Algorithm \ref{alg:SF_algo}) at off-grid points $Y=\{\mathbf{y}_i\}_{i=1}^{N}$ that do not coincide
with the uniform Cartesian grid points $X=\{\mathbf{x}_i\}_{i=1}^{N}$ of the domain. KSF-Chol replaces the original linear
interpolation of these off-grid evaluations with kernel interpolation. We consider a scalar velocity component
$u:\Omega\rightarrow\mathbb{R}$ of the respective velocity vector field $\mathbf{u}=(u,v,w)$ sampled at the $N$ observation points
$
X
=
\{\mathbf{x}_1,\ldots,\mathbf{x}_N\}
\subset\Omega 
$ where 
$
\mathbf{x}_i\in\mathbb{R}^{d},
$ and $ i=1,\ldots,N,
$
with corresponding sampled values
$
\mathbf{u}_X
=
\begin{pmatrix}
u(\mathbf{x}_1) & \cdots & u(\mathbf{x}_N)
\end{pmatrix}^{T}
\in\mathbb{R}^{N}.
$
Thus, each $\mathbf{x}_i$ denotes the $d$-dimensional coordinate vector of
one observation point. For a positive-definite kernel
$K:\Omega\times\Omega\rightarrow\mathbb{R}$, the representer theorem
\cite{ScholkopfHerbrichSmola2001} implies that the corresponding kernel
interpolant admits the representation
\begin{equation}
\psi(\mathbf{u}_X)(\mathbf{y})
=
\sum_{j=1}^{N}
\alpha_j K(\mathbf{y},\mathbf{x}_j)
=
K(\mathbf{y},X)K(X,X)^{-1}\mathbf{u}_X,
\qquad
\mathbf{y}\in\Omega,
\label{eq:closed-form-kernel-interpolation}
\end{equation}
where $\mathbf{y}\in\mathbb{R}^{d}$ denotes an arbitrary evaluation point and
\begin{equation}
K(\mathbf{y},X)
=
\begin{pmatrix}
K(\mathbf{y},\mathbf{x}_1) &
\cdots &
K(\mathbf{y},\mathbf{x}_N)
\end{pmatrix}
\in\mathbb{R}^{1\times N}
\end{equation}
is the corresponding kernel row vector. Furthermore,
$
K(X,X)_{ij}
=
K(\mathbf{x}_i,\mathbf{x}_j),
$ with $
i,j=1,\ldots,N.
$
Thus, $K(X,X)\in\mathbb{R}^{N\times N}$ is the kernel matrix associated
with the observation points (\textit{kernel centers}), which are in our specific case equidistant Cartesian grid points, and
the kernel coefficients
$
\boldsymbol{\alpha}
=
\begin{pmatrix}
\alpha_1 & \cdots & \alpha_N
\end{pmatrix}^{T}
\in\mathbb{R}^{N}
$
are obtained from the \textit{kernel system}
\begin{equation}
K(X,X)\boldsymbol{\alpha}
=
\mathbf{u}_X.
\label{eq:stable-fluids-kernel-system}
\end{equation}
At every advection step, one departure point is associated with each of the
$N$ cell centers. Hence, we denote the corresponding set of $N$
evaluation points by
$
Y
=
\{\mathbf{y}_1,\ldots,\mathbf{y}_N\}
\subset\Omega,
$ with $
\mathbf{y}_i\in\mathbb{R}^{d},
$
where the off-grid locations $\mathbf{y}_i$ generally do not coincide with the
Cartesian observation points $\mathbf{x}_i$. Once the coefficients are known, all
velocities at the evaluation points are obtained simultaneously by the \textit{kernel evaluation}
\begin{equation}
\mathbf{u}(Y)
=
K(Y,X)\boldsymbol{\alpha},
\qquad
K(Y,X)\in\mathbb{R}^{N\times N},
\label{eq:kernel-evaluation-query-points}
\end{equation}
with
$
K(Y,X)_{ij}
=
K(\mathbf{y}_i,\mathbf{x}_j).
$
For the semi-Lagrangian advection step, the observation points are the fixed cell centers, denoted by $X_{\mathrm{grid}}$, while
$Y_{\mathrm{mid}}$ and $Y_{\mathrm{dep}}$ denote the intermediate and final
departure points. Using a second-order Runge–Kutta scheme 
(RK2) for backtracking, they are computed by
\begin{equation}
Y_{\mathrm{mid}}
=
X_{\mathrm{grid}}
-
\frac{\Delta t}{2}\mathbf{u}(X_{\mathrm{grid}}),
\qquad
Y_{\mathrm{dep}}
=
X_{\mathrm{grid}}
-
\Delta t\,\mathbf{u}(Y_{\mathrm{mid}}).
\label{eq:semi-lagrangian-departure-points}
\end{equation}
Here, $\mathbf{u}(X_{\mathrm{grid}})$ denotes the velocity values that are already
available at the cell centers. The RK2 scheme uses these values to
construct the intermediate points $Y_{\mathrm{mid}}$. Then the
off-grid velocity $\mathbf{u}(Y_{\mathrm{mid}})$ is evaluated by the
kernel interpolation. A second interpolation $\mathbf{u}(Y_{\mathrm{dep}})$ gives the velocity at the desired backtracked points.

\subsection{Sparse Cholesky Preconditioner}
\label{sec:sparse-cholesky-preconditioning}

The sparse Cholesky construction from \cite{sparse_chol_leibler} provides a reverse-maximin permutation
$\mathbf{P}$ and a sparse inverse Cholesky factor $\mathbf{L}$ such that an inverse kernel matrix can be approximated by 
$
\mathbf{K}^{-1}
    \approx
    \mathbf{P}^{T}
    \mathbf{L}\mathbf{L}^{T}
    \mathbf{P}
    \label{eq:inverse-cholesky-original-ordering}
$.
The sparsity parameter $\rho_{\mathrm{chol}}$ determines the accuracy of this approximation.
The preconditioner used in the PCG iteration is consequently defined as
\begin{equation}
    \mathbf{M}^{-1}
    =
    \mathbf{P}^{T}
    \mathbf{L}\mathbf{L}^{T}
    \mathbf{P}.
    \label{eq:sparse-inverse-cholesky-preconditioner}
\end{equation}
It is not assembled explicitly. Instead, its application to a residual vector
$\mathbf{r}$ is evaluated in factored form as
$
    \mathbf{z}
    =
    \mathbf{M}^{-1}\mathbf{r}
    =
    \mathbf{P}^{T}
    \left(\mathbf{L}
    \left(
        \mathbf{L}^{T}
        (\mathbf{P}\mathbf{r})
    \right)\right).
    \label{eq:preconditioner-application}
$
One application requires two sparse matrix--vector products and two
permutations with an overall complexity of $\mathcal{O}(N\rho_{\mathrm{chol}}^d)$. The factor and the permutation are built once and reused in every time step. The preconditioner application
is implemented with \texttt{CuPy} and the sparse factor $\mathbf{L}$ is cached in compressed sparse row (CSR) format with \texttt{scipy.sparse}.

\paragraph{FFT based application of the Kernel Matrix}
\label{sec:fft-matern-kernel-application}

In the case of the uniform Cartesian grid, the kernel matrix entries of the stationary and isotropic Gaussian and Matérn kernels depend only on the relative grid displacement of the kernel centers, such that
$
K_{\boldsymbol{i},\boldsymbol{j}}
=
g_{\boldsymbol{i}-\boldsymbol{j}},
$
and therefore form a symmetric three-level Toeplitz matrix whose matrix--vector product is a three-dimensional discrete convolution that can be evaluated efficiently by embedding it into a multilevel circulant matrix and applying multidimensional FFTs \cite{ChanJin2007,RakhubaOseledets2015}.

\subsection{Near-field/far-field evaluation of the kernel expansion}
\label{sec:cheb-svd-fft-evaluation}

After the kernel coefficients
$\boldsymbol{\alpha}_{u}$,
$\boldsymbol{\alpha}_{v}$, and
$\boldsymbol{\alpha}_{w}$
have been calculated by PCG, the remaining operation required by the
semi-Lagrangian advection step is the evaluation of the kernel expansion at
off-grid target points $\mathbf{y}_i$. For one velocity component this operation is
\begin{equation}
    f(\mathbf{y}_i)
    =
    \sum_{j=1}^{N}
    K(\mathbf{y}_i-\mathbf{x}_j)\alpha_j,
    \qquad
    i=1,\ldots,N.
    \label{eq:offgrid-kernel-expansion}
\end{equation}
A direct evaluation requires $\mathcal{O}(N^2)$ kernel evaluations. We
instead combine an exact near field with a Chebyshev approximation and a FFT-based convolution in the far field. Our construction builds upon the fast multipol methods and near-field/far-field
splitting introduced by
\cite{GreengardRokhlin1987,FongDarve2009,MessnerBramasCoulaudDarve2012}.

\paragraph{Target decomposition and near--far splitting}
The source points $\mathbf{x}_j$ are the cell centers of the Cartesian
grid. We denote the set of integer lattice indices by
$
    \Omega_h
    =
    \{0,\ldots,N_x-1\}
    \times
    \{0,\ldots,N_y-1\}
    \times
    \{0,\ldots,N_z-1\}
    \subset\mathbb{Z}^3.
    \label{eq:grid-index-set}
$
Equivalently, a source point is written as
$\mathbf{x}_{\mathbf{j}}$, where
$\mathbf{j}=(j_x,j_y,j_z)^T\in\Omega_h$
is its integer lattice index. For each off-grid target $\mathbf{y}_i$, let
$\mathbf{b}_i=(b_{i,x},b_{i,y},b_{i,z})^T\in\Omega_h$
be the integer lattice index of the Cartesian cell containing the target, and
let $\mathbf{x}_{\mathbf{b}_i}$ denote the center of this cell. Then
\begin{equation}
    \mathbf{y}_i
    =
    \mathbf{x}_{\mathbf{b}_i}
    +
    \boldsymbol{\delta}_i,
    \qquad
    |\delta_{i,q}|
    \leq
    \frac{\Delta q}{2},
    \qquad
    q\in\{x,y,z\}.
    \label{eq:target-cell-decomposition}
\end{equation}
Introducing normalized local coordinates
\begin{equation}
    \mathbf{t}_i
    =
    \left(
        \frac{2\delta_{i,x}}{\Delta x},
        \frac{2\delta_{i,y}}{\Delta y},
        \frac{2\delta_{i,z}}{\Delta z}
    \right)^T
    \in[-1,1]^3,
\end{equation}
and the integer lattice offset
$
    \mathbf{m}
    =
    \mathbf{b}_i-\mathbf{j}
    \in\mathbb{Z}^3,
    $ with $
    \mathbf{b}_i,\mathbf{j}\in\Omega_h,
    \label{eq:lattice-offset}
$
the translation invariance of the stationary Mat\'ern kernel gives
\begin{equation}
    K(\mathbf{y}_i-\mathbf{x}_{\mathbf{j}})
    =
    K\left(
        \mathbf{r}_{\mathbf{m}}
        +
        \mathbf{H}\mathbf{t}_i
    \right),
    \label{eq:lag-local-kernel}
\end{equation}
where the displacement associated with the lattice-offset index
$\mathbf{m}$ is
\begin{equation}
    \mathbf{r}_{\mathbf{m}}
    =
    \left(
        m_x\Delta x,
        m_y\Delta y,
        m_z\Delta z
    \right)^T,
    \end{equation}
 and $
    \mathbf{H}
    =
    \frac{1}{2}
    \operatorname{diag}
    (\Delta x,\Delta y,\Delta z).
$
The set of all lattice offsets that can occur on the finite Cartesian grid is
$
    \mathcal{M}_h
    =
    \left\{
        \mathbf{b}-\mathbf{j}:
        \mathbf{b},\mathbf{j}\in\Omega_h
    \right\}.
    \label{eq:all-lattice-offsets}
$
The off-grid dependence is restricted to the cube
$[-1,1]^3$, while the grid dependence is parameterized by the
discrete lattice-offset index $\mathbf{m}\in\mathcal{M}_h$.
This is the central property exploited by the following approach.

For an integer $c\geq0$, define the near-field set
$
    \mathcal{N}_{c}
    =
    \left\{
        \mathbf{m}\in\mathcal{M}_h:
        \|\mathbf{m}\|_{\infty}\leq c
    \right\},
    \label{eq:near-lag-set}
$
and the complementary far-field offset set
$
    \mathcal{F}_{c}
    =
    \mathcal{M}_h\setminus\mathcal{N}_c.
    \label{eq:far-offset-set}
$
The $(2c+1)^3$ near interactions are evaluated directly with the original
kernel and only the complementary far field
$\mathbf{m}\in\mathcal{F}_c$ is approximated. This avoids
approximating the rapidly varying kernel close to the origin
\cite{GreengardRokhlin1987,FongDarve2009}.

\paragraph{Chebyshev representation of the local target displacement}
For each far-field lattice offset
$\mathbf{m}\in\mathcal{F}_c$, the remaining dependence on
$\mathbf{t}\in[-1,1]^3$ is approximated by Chebyshev polynomials
\cite{Trefethen2013,FongDarve2009}. Let
$p_{\mathrm{cheb}}\in\mathbb{N}$ denote the maximum total polynomial degree,
referred to below as the Chebyshev order. We use the total-degree index set
$
    \mathcal{I}_{p_{\mathrm{cheb}}}
    =
    \left\{
        \boldsymbol{\nu}\in\mathbb{N}_0^3:
        |\boldsymbol{\nu}|_1\leq p_{\mathrm{cheb}}
    \right\},
$
with
$
    P
    =
    |\mathcal{I}_{p_{\mathrm{cheb}}}|
    =
    \binom{p_{\mathrm{cheb}}+3}{3},
    \label{eq:total-degree-size}
$
where $P$ is the number of retained multi-indices and therefore the number
of three-dimensional Chebyshev basis functions. Enumerating
$\mathcal{I}_{p_{\mathrm{cheb}}}=\{\boldsymbol{\nu}^{(q)}\}_{q=1}^{P}$, the corresponding
basis functions are
$
    \Psi_{\boldsymbol{\nu}}(\mathbf{t})
    =
    T_{\nu_x}(t_x)
    T_{\nu_y}(t_y)
    T_{\nu_z}(t_z),
$
where $T_k$ denotes the Chebyshev polynomial of the first kind. We collect
these functions in
$
    \boldsymbol{\Psi}(\mathbf{t})
    =
    \left(
        \Psi_{\boldsymbol{\nu}^{(1)}}(\mathbf{t}),
        \ldots,
        \Psi_{\boldsymbol{\nu}^{(P)}}(\mathbf{t})
    \right)^T
    \in\mathbb{R}^{P}.
    \label{eq:cheb-basis-vector}
$
Hence, for every far-field lattice offset
$\mathbf{m}\in\mathcal{F}_c$,
\begin{equation}
    K\left(
        \mathbf{r}_{\mathbf{m}}
        +
        \mathbf{H}\mathbf{t}
    \right)
    \approx
    \sum_{\boldsymbol{\nu}\in\mathcal{I}_{p_{\mathrm{cheb}}}}
    C_{\boldsymbol{\nu},\mathbf{m}}
    \Psi_{\boldsymbol{\nu}}(\mathbf{t})
    =
    \boldsymbol{\Psi}(\mathbf{t})^T
    \mathbf{c}_{\mathbf{m}},
    \label{eq:cheb-far-expansion}
\end{equation}
where
$\mathbf{c}_{\mathbf{m}}\in\mathbb{R}^{P}$
is the vector of total-degree Chebyshev coefficients
$C_{\boldsymbol{\nu},\mathbf{m}}$ for the lattice offset $\mathbf{m}$.

\paragraph{Chebyshev coefficient fit}
The interpolation matrix required to determine these coefficients is
constructed once before time integration. In each coordinate direction we
use the $p_{\mathrm{cheb}}+1$ Chebyshev--Lobatto nodes
$
    t_k
    =
    \cos\left(\frac{k\pi}{p_{\mathrm{cheb}}}\right),
$ with $
    k=0,\ldots,p_{\mathrm{cheb}}.
    \label{eq:app-cheb-lobatto-1d}
$
Their three-dimensional tensor product contains
$
    L=(p_{\mathrm{cheb}}+1)^3
$
nodes, which are enumerated as
$\{\mathbf{t}^{(\ell)}\}_{\ell=1}^{L}\subset[-1,1]^3$.
The overdetermined Chebyshev interpolation matrix
$\mathbf{D}\in\mathbb{R}^{L\times P}$ is defined by
\begin{equation}
    D_{\ell q}
    =
    \Psi_{\boldsymbol{\nu}^{(q)}}
    \left(\mathbf{t}^{(\ell)}\right),
    \qquad
    \ell=1,\ldots,L,
    \qquad
    q=1,\ldots,P.
    \label{eq:app-cheb-interpolation-matrix}
\end{equation}
For a far-field lattice offset
$\mathbf{m}\in\mathcal{F}_c$, define the vector of kernel samples
$\mathbf{f}_{\mathbf{m}}\in\mathbb{R}^{L}$ by
\begin{equation}
    \left[\mathbf{f}_{\mathbf{m}}\right]_{\ell}
    =
    K\left(
        \mathbf{r}_{\mathbf{m}}
        +
        \mathbf{H}\mathbf{t}^{(\ell)}
    \right),
    \qquad
    \ell=1,\ldots,L.
    \label{eq:app-cheb-kernel-lobatto-samples}
\end{equation}
The associated coefficient vector
$\mathbf{c}_{\mathbf{m}}\in\mathbb{R}^{P}$ is obtained from the
least-squares problem
\begin{equation}
    \mathbf{c}_{\mathbf{m}}
    =
    \underset{\mathbf{c}\in\mathbb{R}^{P}}{\operatorname{argmin}}
    \left\|
        \mathbf{D}\mathbf{c}
        -
        \mathbf{f}_{\mathbf{m}}
    \right\|_2
    =
    \mathbf{D}^{\dagger}
    \mathbf{f}_{\mathbf{m}},
    \label{eq:app-cheb-coefficient-fit}
\end{equation}
where
$\mathbf{D}^{\dagger}\in\mathbb{R}^{P\times L}$
is the Moore--Penrose pseudoinverse of the interpolation matrix
$\mathbf{D}$.

\paragraph{Sampling and SVD compression}
Forming and storing the coefficient vectors
$\mathbf{c}_{\mathbf{m}}$ for every
$\mathbf{m}\in\mathcal{F}_c$ would require excessive memory. We therefore
calculate a subsampled coefficient matrix only for a representative subset
\begin{equation}
    \mathcal{S}
    =
    \left\{
        \mathbf{m}^{(1)},\ldots,\mathbf{m}^{(S)}
    \right\}
    \subset\mathcal{F}_c,
    \qquad
    |\mathcal{S}|=S,
    \label{eq:app-cheb-far-offset-sample-set}
\end{equation}
where $S$ denotes the number of sampled far-field lattice offsets. All admissible
offsets immediately outside the exact near field,
\begin{equation}
    \mathcal{S}_{\mathrm{shell}}
    =
    \left\{
        \mathbf{m}\in\mathcal{F}_c:
        c<\|\mathbf{m}\|_{\infty}\leq c+4
    \right\},
    \label{eq:app-cheb-deterministic-far-shell}
\end{equation}
are included deterministically. Starting from
$\mathcal{S}_{\mathrm{shell}}$, the remaining offsets are selected
uniformly without replacement from
$\mathcal{F}_c\setminus\mathcal{S}_{\mathrm{shell}}$
until the prescribed sample size $S$ is reached.

Using the kernel-sample vectors
$\mathbf{f}_{\mathbf{m}^{(s)}}$
defined in Eq.~\eqref{eq:app-cheb-kernel-lobatto-samples}, we form the sampled
kernel-value matrix
\begin{equation}
    \mathbf{F}_{\mathcal{S}}
    =
    \left[
        \mathbf{f}_{\mathbf{m}^{(1)}}
        \ \cdots\
        \mathbf{f}_{\mathbf{m}^{(S)}}
    \right]
    \in\mathbb{R}^{L\times S}.
    \label{eq:app-cheb-sampled-kernel-matrix}
\end{equation}
Equivalently, its entries are
\begin{equation}
    \left[\mathbf{F}_{\mathcal{S}}\right]_{\ell s}
    =
    K\left(
        \mathbf{r}_{\mathbf{m}^{(s)}}
        +
        \mathbf{H}\mathbf{t}^{(\ell)}
    \right).
    \label{eq:app-cheb-sampled-kernel-matrix-entry}
\end{equation}
Applying the same pseudoinverse
$\mathbf{D}^{\dagger}$ to all sampled columns gives the sampled
Chebyshev coefficient matrix
\begin{equation}
    \mathbf{C}_{\mathcal{S}}
    =
    \mathbf{D}^{\dagger}
    \mathbf{F}_{\mathcal{S}}
    =
    \left[
        \mathbf{c}_{\mathbf{m}^{(1)}}
        \ \cdots\
        \mathbf{c}_{\mathbf{m}^{(S)}}
    \right]
    \in\mathbb{R}^{P\times S}.
    \label{eq:app-cheb-sampled-cheb-matrix}
\end{equation}

A subsequent singular value decomposition of the sampled coefficient matrix
$
    \mathbf{C}_{\mathcal{S}}
    =
    \mathbf{U}
    \boldsymbol{\Sigma}
    \mathbf{V}^{T}
$ leads to further compression.
Truncating the SVD provides the optimal rank-$R$ approximation in the
Euclidean/Frobenius norm
\cite{EckartYoung1936},
\begin{equation}
    \mathbf{C}_{\mathcal{S}}
    \approx
    \mathbf{U}_{R}
    \boldsymbol{\Sigma}_{R}
    \mathbf{V}_{R}^{T}.
    \label{eq:app-cheb-truncated-sampled-svd}
\end{equation}
The columns of $\mathbf{U}_R\in\mathbb{R}^{P\times R}$ form an orthonormal basis of an $R$-dimensional subspace of the $P$-dimensional Chebyshev coefficient space. The vector $\boldsymbol{\beta}_{\mathbf{m}}\in\mathbb{R}^R$ contains the coordinates of the coefficient vector in this reduced basis, such that $\mathbf{U}_R$ reconstructs the approximate full coefficient vector from only $R$ reduced coefficients
\begin{equation}
    \mathbf{c}_{\mathbf{m}}
    \approx
    \mathbf{U}_R\boldsymbol{\beta}_{\mathbf{m}},
    \qquad
    \mathbf{U}_R:
    \mathbb{R}^R\rightarrow\mathbb{R}^P.
\end{equation}
Unless a fixed rank is prescribed, the implementation selects
\begin{equation}
    R
    =
    \min
    \left(
        R_{\max},
        \#
        \left\{
            r:
            \frac{\sigma_r}{\sigma_1}
            >
            \varepsilon_{\mathrm{svd}}
        \right\}
    \right).
    \label{eq:app-cheb-svd-rank-selection}
\end{equation}

\paragraph{Selected skeleton coordinates and spatial modes}

Since the kernel samples at all $L$ Lobatto nodes satisfy $\mathbf{f}_{\mathbf{m}}\approx\mathbf{D}\mathbf{c}_{\mathbf{m}}$, substituting the reduced representation gives
\begin{equation}
    \mathbf{f}_{\mathbf{m}}
    \approx
    \underbrace{\mathbf{D}\mathbf{U}_R}_{\mathbf{Z}}
    \boldsymbol{\beta}_{\mathbf{m}},
    \qquad
    \mathbf{Z}:
    \mathbb{R}^R\rightarrow\mathbb{R}^L.
\end{equation}
Here, $\mathbf{Z}\in\mathbb{R}^{L\times R}$ maps the compressed coefficients to the kernel values at all $L$ Lobatto nodes. Since $L \gg R$, determining $\boldsymbol{\beta}_{\mathbf{m}}$ from these samples would be unnecessary expensive. To reduce the number of evaluations, a column-pivoted QR applied to $\mathbf{Z}^\top$ selects $R$ numerically independent rows corresponding to $R$ row indices
$s_1,\ldots,s_R$ of $\mathbf{Z}$, and therefore $R$
Chebyshev--Lobatto local coordinates
$
    \left\{
        \mathbf{t}^{(s_r)}
    \right\}_{r=1}^{R}
$ and a square matrix $\mathbf{Z}_s\in\mathbb{R}^{R\times R}$ can be formed. The purpose is to obtain the square evaluation matrix
$
    \mathbf{Z}_{s}
    =
    \left[
        Z_{s_r k}
    \right]_{r,k=1}^{R}
    \in\mathbb{R}^{R\times R}
    \label{eq:app-cheb-skeleton-submatrix}
$
that is nonsingular and numerically suitable for inversion
\cite{BusingerGolub1965}. It suffices to evaluate the kernel only at the corresponding $R$ instead of $L$ skeleton coordinates $\mathbf{B}(\mathbf{m})\in\mathbb{R}^R$. For each selected skeleton local coordinate
$\mathbf{t}^{(s_r)}$, $r=1,\ldots,R$, a complete spatial mode is then
evaluated over every lattice offset
$\mathbf{m}\in\mathcal{M}_h$:
\begin{equation}
    B_r(\mathbf{m})
    =
    \begin{cases}
    K\left(
        \mathbf{r}_{\mathbf{m}}
        +
        \mathbf{H}\mathbf{t}^{(s_r)}
    \right),
    &
    \mathbf{m}\in\mathcal{F}_c,
    \\[0.3em]
    0,
    &
    \mathbf{m}\in\mathcal{N}_c.
    \end{cases}
    \label{eq:app-cheb-skeleton-spatial-mode}
\end{equation}
For a fixed lattice offset $\mathbf{m}$, collect the $R$ skeleton values in
\begin{equation}
    \mathbf{B}(\mathbf{m})
    =
    \left(
        B_1(\mathbf{m}),
        \ldots,
        B_R(\mathbf{m})
    \right)^T
    \in\mathbb{R}^{R}.
    \label{eq:app-cheb-skeleton-value-vector}
\end{equation}
We obtain
\begin{equation}
    \mathbf{B}(\mathbf{m})
    \approx
    \mathbf{Z}_s\boldsymbol{\beta}_{\mathbf{m}}
    \quad\Longrightarrow\quad
    \boldsymbol{\beta}_{\mathbf{m}}
    \approx
    \mathbf{Z}_s^{-1}\mathbf{B}(\mathbf{m}).
\end{equation}

So $\mathbf{Z}_s^{-1}:\mathbb{R}^R\rightarrow\mathbb{R}^R$ maps the reduced number of kernel values to the reduced coefficient coordinates, compressing the required kernel evaluations per far-field offset from $L=(p_{\mathrm{cheb}}+1)^3$ to $R$. Finally, the reduced coefficients can be mapped back to the original coefficient space:
\begin{equation}
    \mathbf{c}_{\mathbf{m}}
    \approx
    \underbrace{\mathbf{U}_R\mathbf{Z}_s^{-1}}_{\mathbf{A}_c}
    \mathbf{B}(\mathbf{m}),
    \qquad
    \mathbf{A}_c:
    \mathbb{R}^R\rightarrow\mathbb{R}^P.
\end{equation}
Inserting this representation into the Chebyshev expansion yields the separated far-field approximation
$
    K(\mathbf{r}_{\mathbf{m}}+\mathbf{H}\mathbf{t})
    \approx
    \boldsymbol{\Psi}(\mathbf{t})^\top
    \mathbf{A}_c\mathbf{B}(\mathbf{m})
$
and motivates the local interpolation weights
$\mathbf{a}(\mathbf{t})\in\mathbb{R}^{R}$ defined by
$
    \mathbf{a}(\mathbf{t})^{T}
    =
    \boldsymbol{\Psi}(\mathbf{t})^{T}
    \mathbf{A}_{\mathrm{c}}
    =
    \boldsymbol{\Psi}(\mathbf{t})^{T}
    \mathbf{U}_{R}\mathbf{Z}_s^{-1}.
    \label{eq:local-skeleton-weights}
$
For completeness, define the far-field kernel
$K_{\mathrm{far}}$ on the lattice-offset representation by
\begin{equation}
    K_{\mathrm{far}}
    \left(
        \mathbf{r}_{\mathbf{m}}
        +
        \mathbf{H}\mathbf{t}
    \right)
    =
    \begin{cases}
    K\left(
        \mathbf{r}_{\mathbf{m}}
        +
        \mathbf{H}\mathbf{t}
    \right),
    &
    \mathbf{m}\in\mathcal{F}_c,
    \\[0.3em]
    0,
    &
    \mathbf{m}\in\mathcal{N}_c.
    \end{cases}
    \label{eq:far-kernel-definition}
\end{equation}
One can then write the separated approximation
\begin{equation}
    K_{\mathrm{far}}
    \left(
        \mathbf{r}_{\mathbf{m}}
        +
        \mathbf{H}\mathbf{t}
    \right)
    \approx
    \sum_{r=1}^{R}
    a_r(\mathbf{t})B_r(\mathbf{m}).
    \label{eq:separated-far-kernel}
\end{equation}
The purpose of Eq.~\eqref{eq:separated-far-kernel} is to isolate
the continuous target dependence in $a_r(\mathbf{t})$ from the discrete
lattice-offset dependence in $B_r(\mathbf{m})$. This
separation makes the remaining dependence on
$\mathbf{m}$ amenable to Toeplitz convolution.
\paragraph{FFT evaluation on the Cartesian grid}
For a target $\mathbf{y}_i$ with base-cell index $\mathbf{b}_i$ and local
coordinate $\mathbf{t}_i$, insertion of
Eq.~\eqref{eq:separated-far-kernel} into the far-field part of the kernel
expansion gives
\begin{equation}
\begin{split}
    f_{\mathrm{far}}(\mathbf{y}_i)
    &=
    \sum_{\mathbf{j}\in\Omega_h}
    K_{\mathrm{far}}
    \left(
        \mathbf{r}_{\mathbf{b}_i-\mathbf{j}}
        +
        \mathbf{H}\mathbf{t}_i
    \right)
    \alpha_{\mathbf{j}}
    \\
    &\approx
    \sum_{r=1}^{R}
    a_r(\mathbf{t}_i)
    \sum_{\mathbf{j}\in\Omega_h}
    B_r(\mathbf{b}_i-\mathbf{j})
    \alpha_{\mathbf{j}}.
\end{split}
\label{eq:far-kernel-to-convolution}
\end{equation}
This equation shows why the separated form is introduced. For a given
target, the local coordinate $\mathbf{t}_i$ is inserted only into the
$R$ scalar weights $a_r(\mathbf{t}_i)$ and the complete dependence
on source and target grid indices occurs through the difference
$\mathbf{b}_i-\mathbf{j}$ in $B_r$.
For each spatial mode $B_r$ define
\begin{equation}
    g_r^{(q)}(\mathbf{b})
    =
    \sum_{\mathbf{j}\in\Omega_h}
    B_r(\mathbf{b}-\mathbf{j})
    \alpha_{\mathbf{j}}^{(q)},
    \qquad
    q\in\{u,v,w\}.
    \label{eq:mode-convolution}
\end{equation}
Because $B_r(\mathbf{b}-\mathbf{j})$ depends only on the lattice offset
$\mathbf{b}-\mathbf{j}$ and the kernel centers form a fixed Cartesian
grid, each $B_r$ defines a three-level Toeplitz convolution
\cite{ChanJin2007,RakhubaOseledets2015}. The Fourier representation
$\widehat{B}_r$ of every spatial mode $B_r$ is computed once
before time integration. After zero padding to the size required for a
linear convolution, the mode field
$g_r^{(q)}$ is evaluated by
\begin{equation}
    g_r^{(q)}
    =
    \operatorname{IFFT}_3
    \left(
        \widehat{B}_r
        \odot
        \widehat{\alpha}^{(q)}
    \right),
    \qquad
    r=1,\ldots,R,
    \label{eq:far-field-fft-modes}
\end{equation}
where
$\widehat{\alpha}^{(q)}$
is the three-dimensional FFT of the zero-padded coefficient field
$\alpha^{(q)}$ and $\odot$ denotes pointwise multiplication in Fourier
space.
The same mode fields
$g_r^{(q)}$
are reused for both the Runge--Kutta midpoint scheme and the final off-grid
evaluation. For one velocity component and a target $\mathbf{y}_i$, the
complete kernel expansion is
\begin{equation}
\begin{split}
    f(\mathbf{y}_i)
    \approx{}&
    \sum_{r=1}^{R}
    a_r(\mathbf{t}_i)
    g_r(\mathbf{b}_i)
    \\
    &+
    \sum_{\substack{
        \mathbf{m}\in\mathcal{N}_c\\
        \mathbf{b}_i-\mathbf{m}\in\Omega_h
    }}
    K\left(
        \mathbf{y}_i
        -
        \mathbf{x}_{\mathbf{b}_i-\mathbf{m}}
    \right)
    \alpha_{\mathbf{b}_i-\mathbf{m}} .
    \label{eq:final-cheb-svd-evaluation}
\end{split}
\end{equation}
The first term is the far-field contribution
evaluated by FFT. The second term is the direct, exact evaluation of all kernel interactions
whose lattice offsets satisfy
$\mathbf{m}\in\mathcal{N}_c$.

\paragraph{Parameter selection}
The approximation accuracy is controlled by
$c$, $p_{\mathrm{cheb}}$, $S$, $\varepsilon_{\mathrm{svd}}$, and $R_{\max}$.
The near-field radius $c$ should be increased when the kernel changes rapidly
over one grid cell or when a tighter tolerance is required.

\subsection{KSF-Chol and implementation}
\label{sec:ksf-chol-implementation}

We now summarize the core ideas of KSF-Chol referring to the previously
introduced methods. At the beginning of each time step, the kernel coefficients of the three
velocity components are obtained from
\begin{equation}
    \mathbf{K}\boldsymbol{\alpha}_{u}
    =\mathbf{u}^{\mathrm{force}},
    \qquad
    \mathbf{K}\boldsymbol{\alpha}_{v}
    =\mathbf{v}^{\mathrm{force}},
    \qquad
    \mathbf{K}\boldsymbol{\alpha}_{w}
    =\mathbf{w}^{\mathrm{force}},
\end{equation}
using PCG with the preconditioner
$
    \mathbf{M}^{-1}
    =
    \mathbf{P}^{T}
    \mathbf{L}
    \mathbf{L}^{T}
    \mathbf{P}.
$
Within each PCG iteration, products with $\mathbf{K}$ are applied by
multilevel Toeplitz FFT multiplication.
After the coefficient vectors have been obtained, their Fourier transform is
computed and combined with the $R$ precomputed spatial
modes. The resulting far-field modes are reused for both evaluations
required by the second-order Runge--Kutta backtracking scheme.


\paragraph{PyKeOps parallelization of dense Kernel Evaluation}
A dense \texttt{PyKeOps} implementation serves as an exact
reference for the periodic validation of the near-field/far-field kernel
evaluation.



\paragraph{Computational complexity of KSF-Chol}
Motivated by the scaling $\rho_{\mathrm{chol}}\approx\log(N/\varepsilon_{\mathrm{chol}})$ in \cite{sparse_chol_leibler}, we consider for our asymptotic analysis a sparsity parameter that grows with $N$. For a bounded PCG iteration number $I$ and constant near-field/far-field parameters the computational complexities in time and space are
\begin{equation}
    T_{\mathrm{step}}=\mathcal{O}(N\log^d(N/\varepsilon_{\mathrm{chol}})),
    \qquad
    S_{\mathrm{step}}=\mathcal{O}(N\log^d(N/\varepsilon_{\mathrm{chol}})).
\end{equation}
In the following benchmarks we choose a constant $\rho_{\mathrm{chol}}=7.5$.

\section{Numerical Experiments}
\label{sec:numerical-experiments}
In this section we conduct numerical experiments of KSF-Chol. Section ~\ref{sec:evaluation-methodology} first summarizes the implementation details and
common evaluation criteria.
The individual test cases and their numerical results are then presented in
Section~\ref{sec:test-cases-results}.

\subsection{Implementation and evaluation methodology}
\label{sec:evaluation-methodology}
\subsubsection{Computing Hardware}
All experiments were performed on a workstation equipped with two NVIDIA A40 GPUs, each providing 48\,GB of GPU memory. The system is powered by two AMD EPYC 7352 24-Core processors running at up to 2.30\,GHz, providing a total of 48 physical CPU cores and 48 hardware threads.

\subsubsection{Measures of accuracy}
\label{sec:measures-of-accuracy}
We use the same momentum and divergence residual for all solvers \cite{OberkampfTrucano2002}.
For the two-dimensional case considered here, let $D_x$ and $D_y$ denote
second-order centered finite-difference operators evaluated on the respective grid, and let $\Delta_h$ denote the corresponding discrete Laplace
operator. The momentum residuals at time $t_n$ are then defined with $(\phi,\xi)\in\{(u,x),(v,y)\}$ as 

\begin{equation}
\begin{aligned}
R_{\phi}^{n}
&=
\frac{\phi^{n+1}-\phi^{n-1}}{2\Delta t_{\mathrm{eval}}}
+
u^{n} D_x \phi^{n}
+
v^{n} D_y \phi^{n}
+
\frac{1}{\rho_f} D_{\xi} p^{n}
-
\nu \Delta_h \phi^{n}
-
f_{\xi}^{n}.
\end{aligned}
\label{eq:momentum-residual}
\end{equation}
The magnitude of the momentum residual is quantified by the discrete spatial
$L^2$ norm
\begin{equation}
\left\|
R_{\mathrm{mom}}^{n}
\right\|_{L_h^2}
=
\left[
\Delta x \Delta y
\sum_{(i,j)\in\Omega_f}
\left(
\left(R_{u,i,j}^{n}\right)^2
+
\left(R_{v,i,j}^{n}\right)^2
\right)
\right]^{1/2},
\label{eq:momentum-residual-l2}
\end{equation}
where $\Omega_f$ denotes the set of grid points included in the
residual evaluation. The incompressibility constraint is assessed independently through the
discrete divergence residual
\begin{equation}
R_{\mathrm{div}}^{n}
=
D_x u^{n}
+
D_y v^{n},
\label{eq:divergence-residual}
\end{equation}
with the corresponding discrete $L^2$ norm
\begin{equation}
\left\|
R_{\mathrm{div}}^{n}
\right\|_{L_h^2}
=
\left[
\Delta x \Delta y
\sum_{(i,j)\in\Omega_f}
\left(
R_{\mathrm{div},i,j}^{n}
\right)^2
\right]^{1/2}.
\label{eq:divergence-residual-l2}
\end{equation}

\subsection{Test cases and results}
\label{sec:test-cases-results}
The first benchmark ~\ref{sec:energy-spectrum} analyzes the 2D energy cascade. The following two test cases, presented in
Sections ~\ref{sec:merging-vortices} and~\ref{sec:karman-vortex-street},
address the low-Reynolds-number regime and focus on numerical diffusion and
the conservation of vortical flow structures. Their choice was is inspired by \cite{PicTurb}. The subsequent benchmark in
~\ref{sec:flow-past-sphere}
considers a more turbulent flow and evaluates the ability of KSF-Chol to
reproduce increasingly complex wake dynamics.
\subsubsection{The energy spectrum}
\label{sec:energy-spectrum}

In this first benchmark, we investigate whether the KSF-Chol algorithm correctly reproduces the characteristic $k^{-3}$ energy spectrum of two-dimensional turbulence. We perform a numerical simulation of the incompressible NS equations in vorticity form on a periodic domain.
\paragraph{Setup}
We consider the vorticity--streamfunction $(\omega - \psi)$ formulation of the incompressible Navier--Stokes equations on the two-dimensional torus $\mathbb{T}^2 = [0,2\pi]^2$:
\begin{equation}
\partial_t \omega(\mathbf{x}, t) + \mathbf{u}(\mathbf{x}, t)\cdot \nabla \omega(\mathbf{x}, t)
= \nu \Delta \omega(\mathbf{x}, t) + f(\mathbf{x}), 
\quad \mathbf{x} \in [0,2\pi]^2, \; t \in [0,T],
\label{eq:vorticity-equation}
\end{equation}
\begin{equation}
\omega = -\Delta \psi, 
\quad 
\mathbf{u} = 
\begin{pmatrix}
\partial_y \psi \\
-\partial_x \psi
\end{pmatrix},
\quad 
\nabla \cdot \mathbf{u} = 0,
\quad 
\omega(\mathbf{x},0) = \omega_0(\mathbf{x}).
\label{eq:vorticity-streamfunction}
\end{equation}
Here, $\mathbf{u}$ denotes the divergence-free velocity field, $\omega=\nabla\times\mathbf{u}$ is the corresponding scalar vorticity, and $\psi$ is the streamfunction obtained from the Poisson equation. The initial vorticity is denoted by $\omega_0$, $\nu>0$ is the viscosity coefficient, and $f$ represents the external forcing term.

\paragraph{Energy spectrum}
To analyze the turbulent energy cascade, the kinetic energy spectrum is
computed from the velocity field using its Fourier representation
\cite{Pope2000}.

\begin{figure}
\centering
\includegraphics[width=0.75\textwidth]{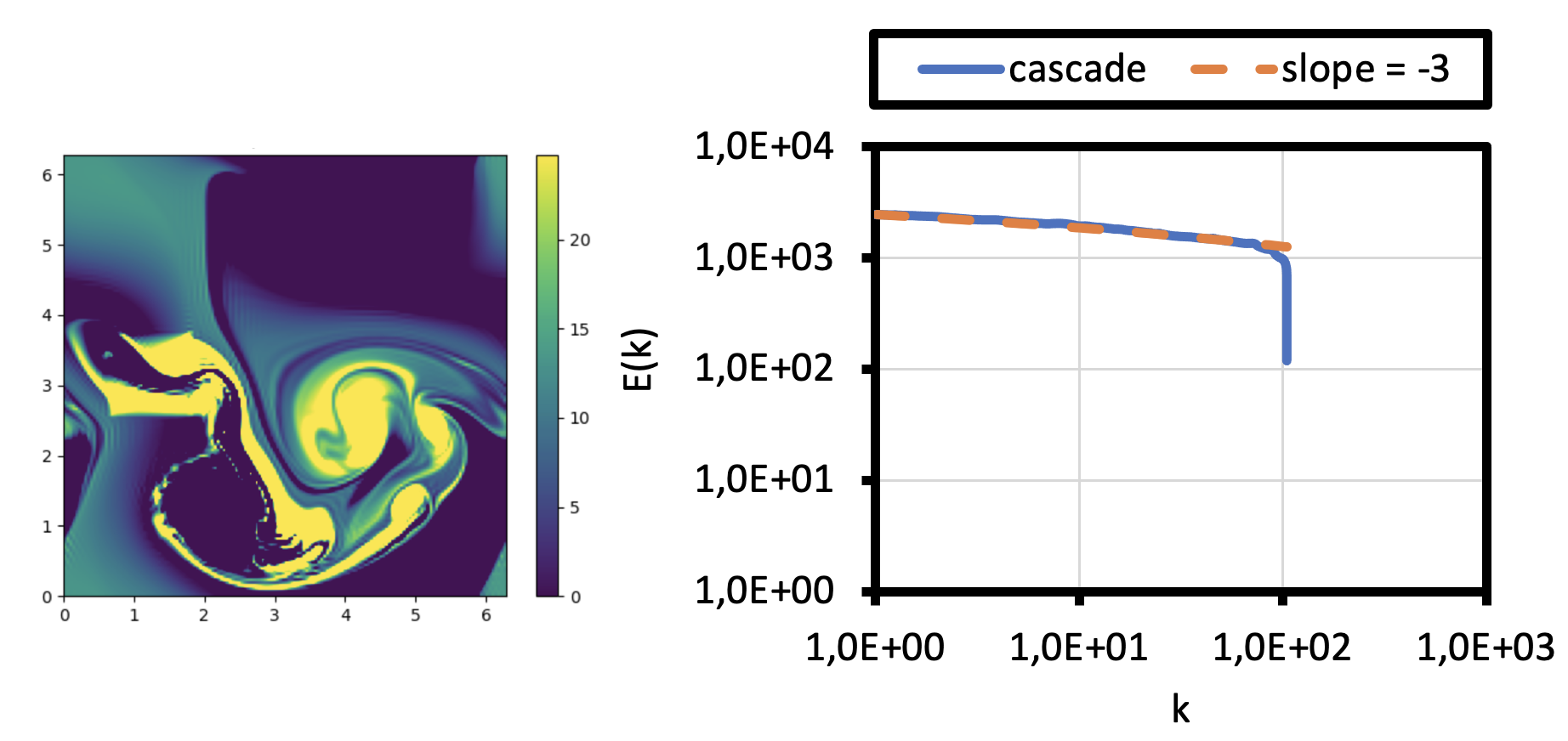}
\caption[]{KSF-Chol solution of the Navier--Stokes equations with local and global forcing. \textbf{Left:} Vorticity $\omega$ at the final time step. \textbf{Right:} The energy spectrum at the final time step and the slope of the expected $k^{-3}$ scaling law are shown.}
\label{fig:energy_casc}
\end{figure}

\paragraph{Results}
According to \autoref{fig:energy_casc}, the energy spectrum approximately follows the $k^{-3}$ scaling law and is therefore consistent with the scaling associated with a two-dimensional forward enstrophy cascade \cite{Kraichnan1967}.

\subsubsection{Two merging vortices}
\label{sec:merging-vortices}

This benchmark investigates how accurately the KSF-Chol algorithm reproduces
the dynamics of two merging vortices. As a numerical reference solution, we
employ a Fourier pseudo-spectral solver. For smooth periodic incompressible
flows, spectral and pseudo-spectral methods provide very high spatial
accuracy, and their accuracy, stability, and convergence properties for
Navier--Stokes problems are well established
\cite{orszag1971accuracy,MadayQuarteroni1982}. Provided that the spatial and
temporal resolutions are sufficiently converged, the resulting solution can
therefore be used as a highly accurate numerical reference.

\paragraph{Error analysis}

\begin{figure}
\centering
\includegraphics[width=0.75\textwidth]{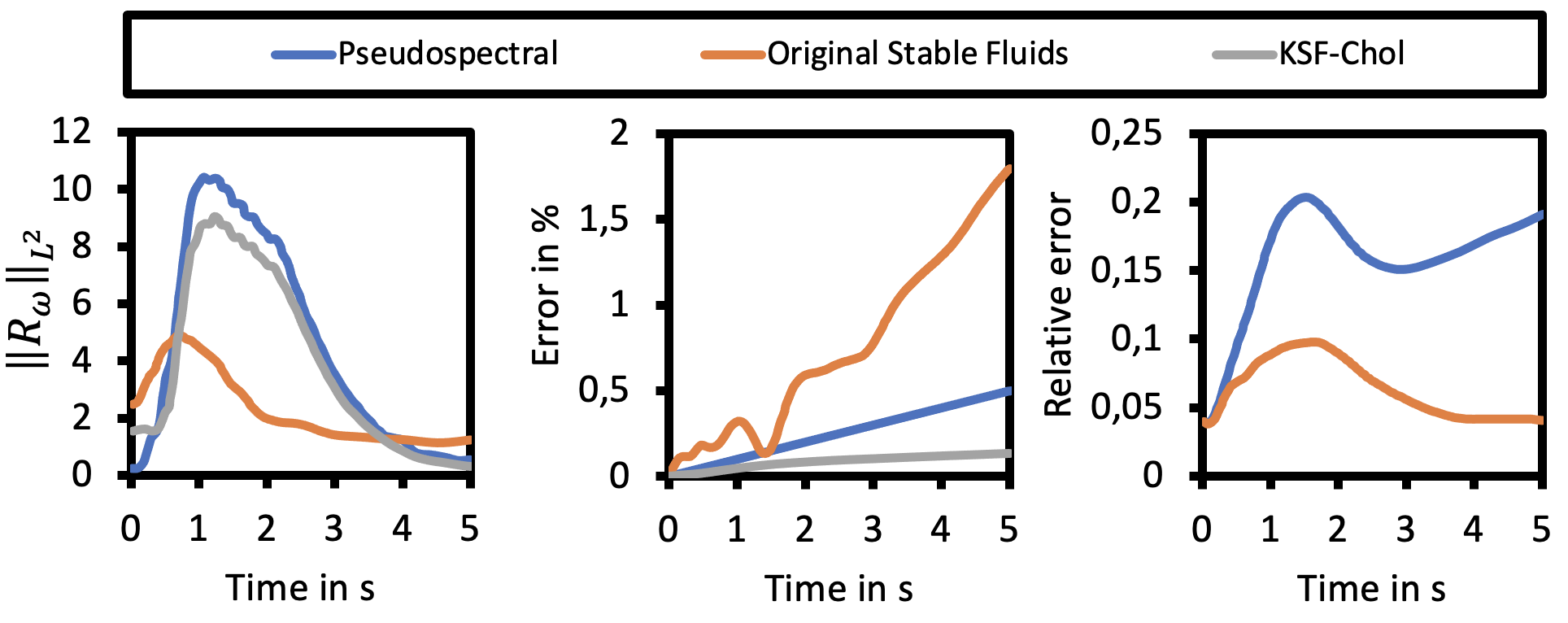}
\caption[]{Error analysis of the two merging vortices simulation with different solver types:  \textbf{Left:} Vorticity-equation residual. \textbf{Middle:} Relative integral-vorticity error. \textbf{Right:} The relative error in the $L^2$ norm of the velocity fields of the KSF-Chol and original Stable Fluids algorithm with respect to the pseudospectral reference solution is plotted for every time step of the simulation. (see \autoref{tab:merging_vort_tab} for setup details) }
\label{fig:merging_vort_diag}
\end{figure}

In this benchmark we analyze the conservation of integral vorticity
and the residual of the vorticity equation.
The integral vorticity is expected to remain
constant in time. It is approximated on the computational grid by
$
M_n
=
\sum_{i,j}
\omega_{i,j}^{n}\,\Delta A
$ with $
\Delta A=\Delta x\Delta y,
\label{eq:total-vorticity}
$. Its relative deviation
from the initial integral vorticity value is defined as
\begin{equation}
\mathrm{err}_n
=
100\,
\frac{|M_n-M_0|}{|M_0|},
\label{eq:total-vorticity-error}
\end{equation}
and shown in \autoref{fig:merging_vort_diag}. This quantity provides a
measure of the numerical loss or gain of integral vorticity.

In addition, the consistency of the numerical solution with the vorticity
equation is measured by
\begin{equation}
R_{\omega}^{n}
=
\frac{
\omega^{n+1}-\omega^{n-1}
}{
2\Delta t_{\mathrm{eval}}
}
+
u^{n}D_x\omega^{n}
+
v^{n}D_y\omega^{n}
-
\nu\Delta_h\omega^{n},
\label{eq:vorticity-residual}
\end{equation}
where $D_x$, $D_y$, and $\Delta_h$ are the same second-order centered
operators used in the residual evaluation described in
Section~\ref{sec:measures-of-accuracy}.

The magnitude of the vorticity residual is quantified by the discrete
spatial $L^2$ norm
\begin{equation}
\left\|R_{\omega}^{n}\right\|_{L_h^2}
=
\left[
\Delta x\,\Delta y
\sum_{i=1}^{N_x}
\sum_{j=1}^{N_y}
\left(R_{\omega,i,j}^{n}\right)^2
\right]^{1/2}.
\label{eq:vorticity-residual-l2}
\end{equation}

\paragraph{Results}

\begin{figure}[t]
\centering
\includegraphics[width=0.95\textwidth]{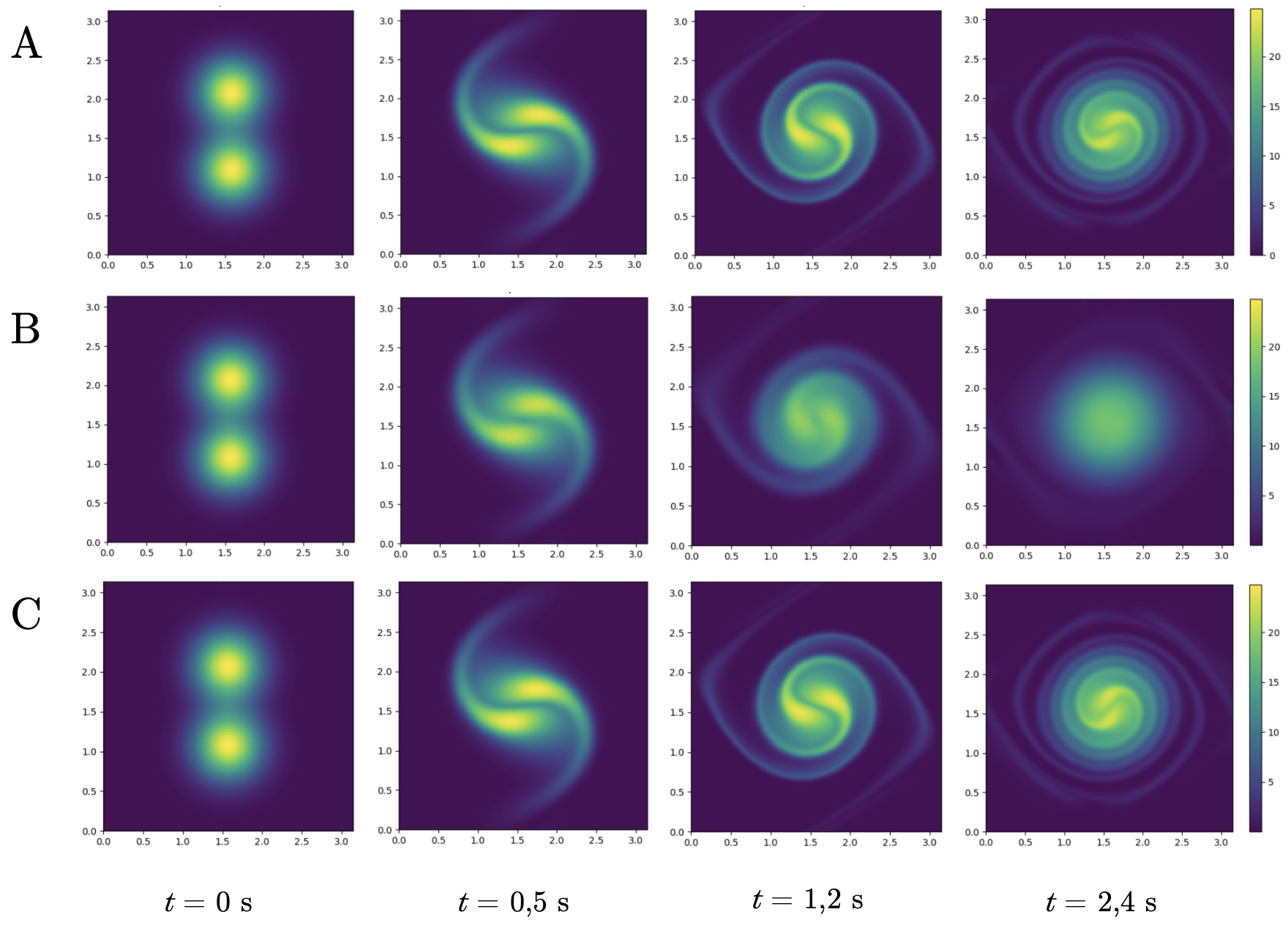}
\caption[]{Plots of the vorticity $\omega$ at different time steps for the two merging vortices benchmark with three solver types. \textbf{A:} Pseudospectral \textbf{B:} Stable Fluids \textbf{C:} KSF-Chol.}
\label{fig:merging_vort_pic}
\end{figure}

The results in \autoref{fig:merging_vort_pic} show that the KSF-Chol method reproduces the reference solution more accurately than the original SF approach. In particular, the relative $L^2$ error remains consistently smaller as shown in \autoref{fig:merging_vort_diag}. The integral-vorticity error in \autoref{fig:merging_vort_diag} is smaller for KSF-Chol than for Stable Fluids, and the vorticity-equation residual is comparable to the pseudo-spectral reference residual. For this benchmark, KSF-Chol remains stable at $\Delta t=0.01$, whereas the pseudo-spectral reference requires $\Delta t=0.001$. The larger stable time step is obtained at the expense of a higher wall-clock cost per time step.

\subsubsection{Von Kármán vortex street}
\label{sec:karman-vortex-street}
In this benchmark, it is investigated how accurately
KSF-Chol reproduces the von K\'arm\'an vortex street behind a cylindrical
obstacle. As a numerical reference, we choose a finite-difference solver
based on the projection/fractional-step framework introduced by
~\cite{Chorin1968}. Projection methods have since been developed into
a broad class of solvers ~\cite{GuermondMinevShen2006}. Modified variants of this
general solver architecture are used for high-fidelity simulations,
including direct numerical simulations of incompressible turbulent
flows~\cite{Costa2018,VanDerPoel2015}.

\paragraph{Setup}

We consider the incompressible Navier--Stokes equations in velocity formulation on a rectangular domain $\Omega$ as stated in ~\eqref{eq:inkomp_eq} but with zero forcing.

\paragraph{Error analysis}

\begin{figure}[t]
\centering
\includegraphics[width=0.95\textwidth]{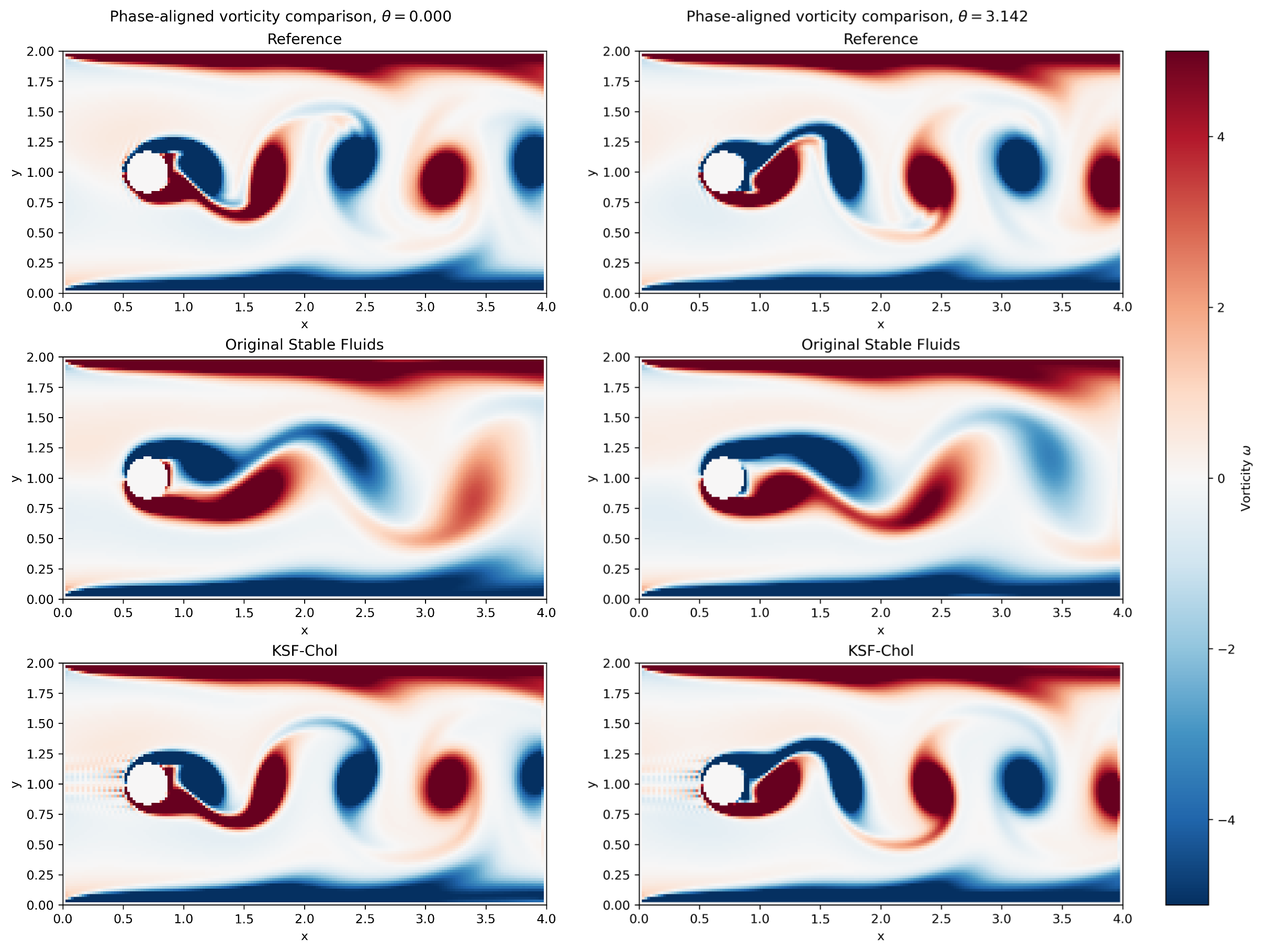}
\caption[]{Phase-aligned vorticity fields of the Reference, Stable Fluids, and KSF-Chol solvers at $\theta=0$ (left) and $\theta=3.142$ (right). The fields are sampled at corresponding normalized phases of their respective vortex-shedding cycles.}
\label{fig:pic_karman}
\end{figure}

\begin{figure}[t]
\centering
\includegraphics[width=0.65\textwidth]{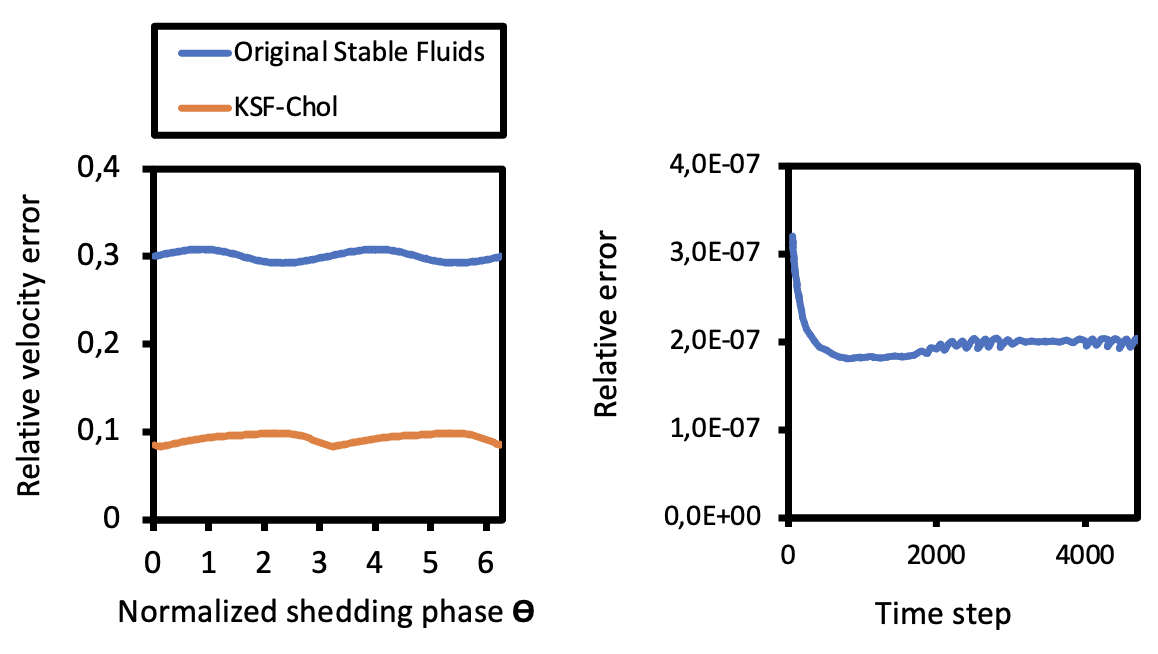}
\caption[]{\textbf{Left:} Phase-conditioned relative $L^2$ velocity error of the Stable Fluids and KSF-Chol solvers with respect to the reference solution as a function of the normalized shedding phase $\theta$. Each solver is sampled at the corresponding phase of its vortex-shedding cycle. \textbf{Right:} Relative $L^2$ error between the advected velocity fields obtained with dense and near-field/far-field kernel evaluation. The error is evaluated after the complete midpoint-RK2 advection step, including both kernel interpolations.}
\label{fig:phase_conditioned_velocity_error_three_solver}
\end{figure}

The numerical solutions are compared using the momentum
and divergence residuals introduced in
Section~\ref{sec:measures-of-accuracy}.
The vortex-shedding frequency and the accuracy of the phase aligned flow
field are assessed separately. The Strouhal number measures whether a solver reproduces
the correct shedding frequency. A direct comparison of the velocity fields at
equal physical times would intertwine errors of the shedding frequency with
differences in the flow field and even a small period error would produce
an increasing phase shift over successive cycles. For the phase-conditioned analysis, one complete vortex shedding cycle
is selected for the reference solution and for each SF
solver. Introducing the normalized shedding phase
$
    \theta \in [0,2\pi],
$ with $
    s=\frac{\theta}{2\pi},
$
the corresponding physical times are
$
    t_{\mathrm{ref}}(\theta)
    =
    t_{\mathrm{ref},0}
    +
    s T_{\mathrm{ref}},
$ and $
    t_{\mathrm{sol}}(\theta)
    =
    t_{\mathrm{sol},0}
    +
    s T_{\mathrm{sol}}.
$
The velocity fields are linearly interpolated at these times. For each phase,
the relative velocity error shown in
Fig.~\ref{fig:phase_conditioned_velocity_error_three_solver} is computed over
the fluid domain $\Omega_f$ as
\begin{equation}
    e_u(\theta)
    =
    \frac{
        \left\|
        \mathbf{u}_{\mathrm{sol}}(\theta)
        -
        \mathbf{u}_{\mathrm{ref}}(\theta)
        \right\|_{L^2(\Omega_f)}
    }{
        \left\|
        \mathbf{u}_{\mathrm{ref}}(\theta)
        \right\|_{L^2(\Omega_f)}
    },
\end{equation}
with the discrete norm being evaluated using the common cell area $\Delta A$. The scalar velocity cycle error reported in the table combines these
phase-resolved errors over the complete normalized cycle and is defined as
\begin{equation}
E_{u,\mathrm{cycle}}
=
\frac{
\left\|
\mathbf{u}_{\mathrm{sol}}
-
\mathbf{u}_{\mathrm{ref}}
\right\|_{L^2([0,2\pi]\times\Omega_f)}
}{
\left\|
\mathbf{u}_{\mathrm{ref}}
\right\|_{L^2([0,2\pi]\times\Omega_f)}
}.
\end{equation}
The corresponding phase-aligned vorticity fields in
Fig.~\ref{fig:pic_karman} are reconstructed from the phase-conditioned
velocity fields.

\paragraph{Dense vs. Near-field/far-field Kernel Evaluation}
The accuracy of the near-field/far-field kernel evaluation is assessed using the dense kernel evaluation as a reference. Every 50 time steps, the interpolation for $\mathbf{U}_{\mathrm{mid}}$ and $\mathbf{U}^{\mathrm{adv}}$ is evaluated independently using the near-field/far-field and dense kernel evaluations. The two resulting advected velocity fields are then compared in $L^2$ norm according to
\begin{equation}
e_{\mathrm{rel},L^2}
=
\frac{
\left\|
\mathbf{u}_{\mathrm{approx}}^{\mathrm{adv}}
-
\mathbf{u}_{\mathrm{dense}}^{\mathrm{adv}}
\right\|_2
}{
\left\|
\mathbf{u}_{\mathrm{dense}}^{\mathrm{adv}}
\right\|_2
},
\end{equation}
As shown in \autoref{fig:phase_conditioned_velocity_error_three_solver}, the relative $L^2$ error remains close to $2\times10^{-7}$. This nearly constant error indicates that the near-/far-field evaluation approximates the dense kernel evaluation with high accuracy in the given vortex-shedding regime. The Strouhal number obtained with the near-/far-field evaluation differs from the dense result by only $10^{-6}\%$, further indicating that the approximation preserves the vortex-shedding dynamics of the dense KSF-Chol simulation.

\paragraph{Results}

All three methods require a transient phase before reaching the periodic vortex shedding regime. The original SF method exhibits excessive numerical diffusion. In contrast, the KSF-Chol method generates a more turbulent flow than the original SF algorithm and reproduces the reference solution more closely as one can observe in \autoref{fig:pic_karman}. Figure~\ref{fig:mom_res_karman} in Appendix~\ref{app:karman-setup} shows that the momentum residual of the reference solver is significantly smaller when compared to both SF solvers. Yet the momentum residual of the KSF-Chol solver is smaller than the one of the original SF algorithm. The divergence residual of both SF methods is nearly equal.

The vortex-shedding statistics in \autoref{tab:karman_parameters} confirm this
difference. Relative to the reference Strouhal number, the original SF method exhibits an error of approximately $19.7\%$, whereas the
KSF-Chol result differs only by $1.3\%$. Hence, KSF-Chol
reproduces the characteristic vortex-shedding frequency substantially more
accurate.

The phase-conditioned velocity comparison in
\autoref{fig:phase_conditioned_velocity_error_three_solver} leads to the
conclusion that over the normalized shedding cycle, the relative $L^2$ velocity
error of KSF-Chol remains below that of the original SF method.
The velocity cycle $L^2$ error reported in
\autoref{tab:karman_parameters} is approximately $9.27\times10^{-2}$ for
KSF-Chol and $3.01\times10^{-1}$ for original SF. Thus, the cycle error of
KSF-Chol is lower
than that of the original SF method. This shows that KSF-Chol not
only reproduces the shedding frequency more accurately, but also yields a
substantially closer agreement with the reference velocity
field throughout the vortex-shedding cycle.

The semi-Lagrangian advection step is not explicitly amenable to the conventional CFL stability restriction. The KSF-Chol implementation nevertheless exhibits practical time-step restrictions at high Reynolds numbers; in our simulations, artifact-free solutions required the CFL number to remain approximately below unity.

\begin{table}[h]
\centering
\caption{Vortex-shedding results and computational performance for the von Kármán vortex street.}
\label{tab:karman_parameters}
\scriptsize
\begin{tabularx}{\linewidth}{@{}l>{\centering\arraybackslash}X>{\centering\arraybackslash}X>{\centering\arraybackslash}X>{\centering\arraybackslash}X@{}}
\toprule
Solver & Reference & Stable Fluids & KSF-Chol & KSF-Chol \\
 & & & (dense) & (approx) \\
\midrule
Memory $L_{\mathrm{sparse}}$ & --- & --- & $0.066\,\mathrm{GB}$ & $0.066\,\mathrm{GB}$ \\
Wall-clock time per time step & $0.000957 \mathrm{s}$ & $0.005908 \mathrm{s}$ & \shortstack{$1.36\,\mathrm{s}$\\with CuPy} & $0.28\,\mathrm{s}$ \\
\midrule
$St$ & $0.24667081$ & $0.19797213$ & $0.24348520$ & $0.24348502$ \\
Velocity cycle $L^2$ error & --- & $3.006190\times10^{-1}$ & $9.269801\times10^{-2}$ & $9.314720\times10^{-2}$ \\
\bottomrule
\end{tabularx}
\end{table}

\subsubsection{3D Flow Past a Sphere at $\mathrm{Re}=1000$}
\label{sec:flow-past-sphere}

We consider a three-dimensional incompressible flow past a stationary sphere at
$\mathrm{Re}=1000$. Flows in this turbulence regime represent a challenging benchmark with a transitional wake that is characterized by an unsteady recirculation region, vortex shedding, and the formation of vortical structures downstream of the sphere. We want to assess if KSF-Chol is able to reproduce the near-wake dynamics as reported by \cite{MarefatAlamPope2022} and ~\cite{RodriguezEtAl2019}.
\cite{MarefatAlamPope2022} investigated the flow past a sphere using large-eddy simulations (LES) combined with an immersed
boundary method (IBM) and a vorticity-based subgrid-scale (SGS) model. Their numerical
results were compared with the direct numerical simulations (DNS) of ~\cite{RodriguezEtAl2019} and experimental data, showing good
agreement for characteristic wake quantities.
In the present benchmark, we focus on one specific quantity reported in the reference
studies. The pressure coefficient $C_{p}$ on the sphere surface,
since it directly reflects the pressure distribution generated by the near-wake
dynamics and flow separation. We furthermore analyze the energy cascade of our simulation.

\paragraph{Setup}
We consider the incompressible Navier--Stokes equations in velocity formulation on a rectangular domain $\Omega$ as stated in ~\eqref{eq:inkomp_eq} but with zero forcing.

\paragraph{KSF-Chol}
For this benchmark only, we apply the Wall-Adapting Local Eddy-Viscosity (WALE) model of ~\cite{NicoudDucros1999} after advection and before molecular diffusion. This benchmark-specific SGS update is not part of the base KSF-Chol algorithm. 

\begin{figure}
\centering
\includegraphics[width=0.65\textwidth]{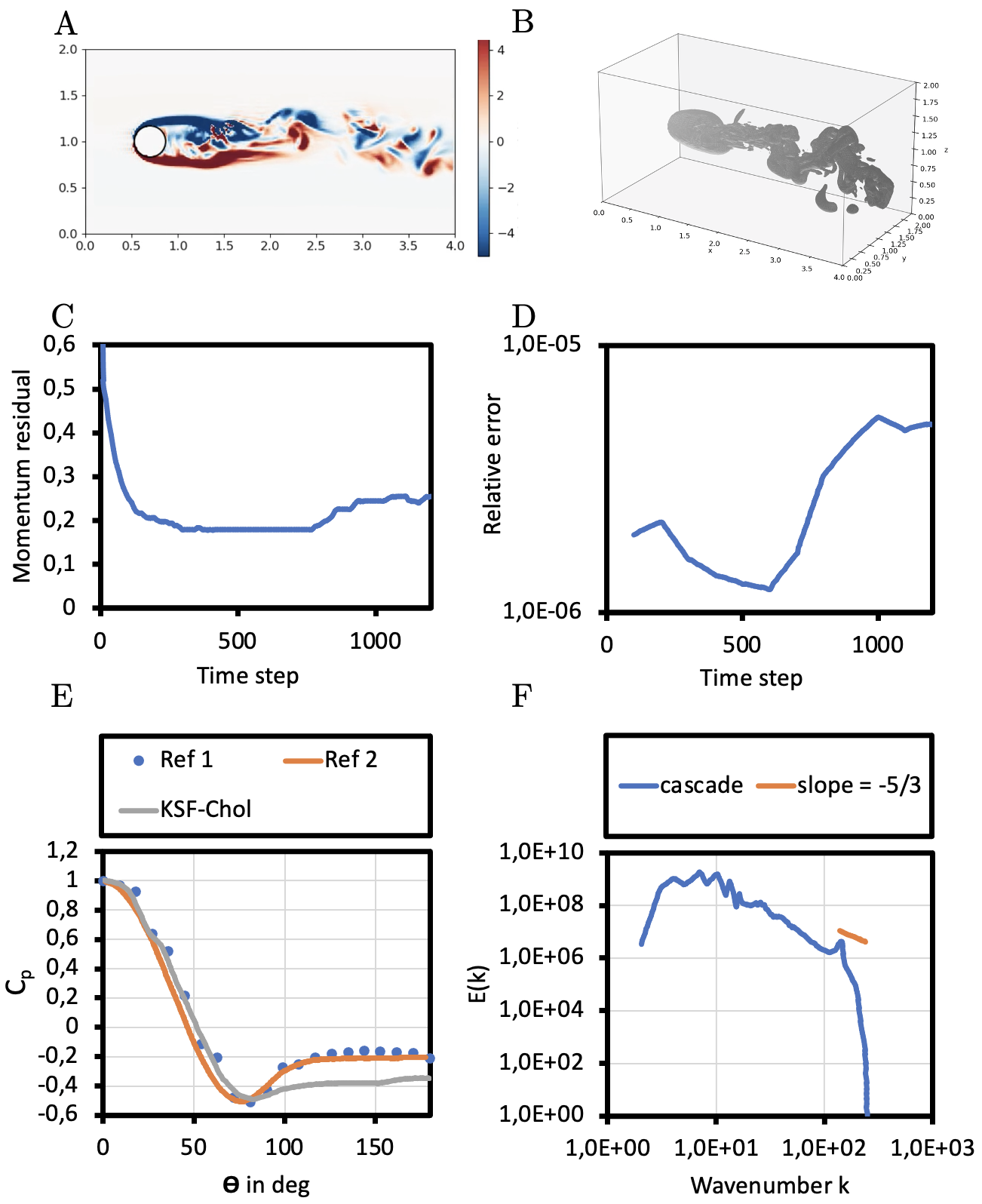}
\caption[]{\textbf{3D flow simulation over a sphere at $\mathrm{Re}=1000$ with KSF-Chol:} \textbf{A:} Midplane vorticity plot at step $1160$. \textbf{B:} 3D vorticity isosurface at step $1160$ \textbf{C:} Momentum residual evolution. \textbf{D:} Relative $L^2$ error between the velocity fields obtained with dense and near-field/far-field kernel evaluation. The error is evaluated after the complete midpoint-RK2 advection step, including both kernel interpolations. \textbf{E:} Distribution of the pressure coefficient $C_{p}$ along the mid-plane section of the sphere surface for KSF-Chol (\textit{grey line}) and the two numerical references using the LE--IBM simulation of ~\cite{MarefatAlamPope2022} (\textit{blue points}) and DNS of ~\cite{RodriguezEtAl2019} (\textit{orange line}). \textbf{F:} The energy spectrum and the slope of the expected $k^{-5/3}$ scaling law are shown.}
\label{fig:re1000_all}
\end{figure}

\paragraph{Results}
The LE--IBM results of ~\cite{MarefatAlamPope2022} are used as the primary numerical reference because their $\mathrm{Re}=1000$ simulation is closest to our setup (regarding the uniform Cartesian domain and the similar grid spacings). ~\cite{MarefatAlamPope2022} also mentions ~\cite{TomboulidesOrszag2000} besides ~\cite{RodriguezEtAl2019} as reference simulation at $\mathrm{Re}=1000$.
The KSF-Chol pressure coefficient shown in \autoref{fig:re1000_all} (\textbf{E}) is compared with the LES--IBM results of ~\cite{MarefatAlamPope2022} and the DNS results of ~\cite{RodriguezEtAl2019}. Following the technique of ~\cite{MarefatAlamPope2022}, the pressure values of KSF-Chol are sampled at different angular positions using points within the radial interval $0.45\leq r/D\leq0.55$ around the sphere surface which lies at $r/D=0.5$. From the front stagnation point to the suction peak, KSF-Chol reproduces the reference $C_{p}$ distribution closely. But downstream of the suction peak the $C_{p}$ distribution of KSF-Chol remains below the reference distribution, corresponding to a weaker pressure recovery and a larger pressure deficit on the rear part of the sphere.
The mean separation angle obtained with KSF-Chol is
$
\theta_{sp}=105.432^\circ .
$
~\cite{MarefatAlamPope2022}, ~\cite{RodriguezEtAl2019}, and ~\cite{TomboulidesOrszag2000} report $102.7^\circ$, $101.4^\circ$, and $102.0^\circ$ at $\mathrm{Re}=1000$. Relative to the primary LES--IBM reference, the deviation of the KSF-Chol result is approximately $2.7\%$. The larger separation angle indicates that KSF-Chol predicts separation slightly further downstream toward the rear stagnation point.

The energy spectrum shown in \autoref{fig:re1000_all} (\textbf{D}) exhibits an approximate $k^{-5/3}$ slope consistent with classical three-dimensional turbulence scaling law \cite{Pope2000}. The residual shown in \autoref{fig:re1000_all}  (\textbf{F}) remains small throughout the simulation while the CFL number is nearly constant, indicating a stable temporal evolution of the KSF-Chol solution.

The computational domain used in this test case consists of
$
N=1{,}458{,}000
$
grid points (\autoref{tab:setup_re1000}). For this simulation, constructing the sparse Cholesky factor requires
$t_{\mathrm{sparse}}=3.1\,\mathrm{h}$ and its storage requires
$13.21\,\mathrm{GB}$ of GPU memory. Thus, the construction and, in particular, the storage of the sparse factor $L_{\mathrm{sparse}}$ constitute a limiting factor for larger domain sizes and would require additional optimization and implementation strategies at substantially higher resolutions. Finally, the substantial acceleration achieved by the near-field/far-field evaluation becomes apparent (details in \autoref{tab:cheb_svd_fft_re1000}), while the approximation reproduces the dense result with a relative $L^2$-error on the order of $10^{-6}$ as one can observe in \autoref{fig:re1000_all} (\textbf{B}).
 A dense evaluation using \texttt{PyKeOps} requires approximately
$14.91\,\mathrm{min}$ per time step, whereas the proposed near-field/far-field approach reduces the wall-clock time to approximately
$10.2\,\mathrm{s}$ per time step.

In particular this 3D test case highlights the significant acceleration of the kernel evaluation achieved by the near-field/far-field approach and demonstrates that, under the assumptions made in \autoref{sec:cheb-svd-fft-evaluation}, large kernel interpolations can be performed rapidly and repeatedly with high accuracy, making the method suitable for the use in arbitrary solver architectures (see Wall-clock time/step in \autoref{tab:performance_re1000}).

\begin{table}[t]
\centering
\caption{Computational performance and separation-angle result for the 3D flow over a sphere benchmark.}
\label{tab:performance_re1000}
\scriptsize
\begin{tabular}{@{}ll@{}}
\toprule
\textbf{Parameter} & \textbf{Value} \\
\midrule
$t_{\mathrm{sparse}}$ (h)           & $3.1\,\mathrm{h}$ \\
Memory $L_{\mathrm{sparse}}$        & $13.21\,\mathrm{GB}$ \\
Wall-clock time/step (approx.)      & $10.2\,\mathrm{s}$ \\
Wall-clock time/step (dense)        & $14.91\,\mathrm{min}$ \\
\midrule
Mean separation angle $\theta_{sp}$ & $105.432^\circ$ \\
\bottomrule
\end{tabular}
\end{table}

\section{Conclusions and Outlook} 
\label{sec:conclusions}
The numerical experiments in this work show that the proposed KSF-Chol algorithm reproduces flow dynamics and turbulence characteristics more accurately than the original SF algorithm and yields results that are closer to those obtained with the reference solvers chosen in this work. This is demonstrated by the analysis of the energy spectra, vortex-shedding frequencies, and relative error of the phase-aligned flow fields of the respective solvers. The numerical diffusion of vorticity produced by KSF-Chol is small compared with accurate reference solvers, addressing a principal limitation of the original SF method. KSF-Chol also provides good estimates of relevant flow quantities, such as the Strouhal number and the mean separation angle, while keeping the stable solver architecture of the original algorithm. Sparse preconditioning and accelerated kernel evaluation give $\mathcal{O}(N\log^d(N/\varepsilon_{\mathrm{chol}}))$ time and space complexity per solver step.

The main finding and advantage of the presented KSF-Chol architecture is the exploitation and simultaneously restriction to a uniform Cartesian domain that allows for fast repeated kernel evaluations with an exact near-field and a compressed far-field with controllable approximation error in every time step of the solver. This idea can therefore be applied to every solver type and is not restricted to our specific fluid application. The KSF-Chol algorithm also allows for extensions, as demonstrated by the incorporation of the SGS model in the 3D sphere benchmark and can accordingly be further modified and improved. Furthermore, the sparse Cholesky factors could be learned within an operator learning framework, eliminating the need to recompute them at every time step. Instead, they can be efficiently predicted in a single forward pass using neural networks or Gaussian processes \cite{Bacho2026Operator}. Since the predicted factors will only be used as preconditioner, the predicion of the factors themselves don't have to be very accurate. 
The restriction of the solver architecture to a uniform Cartesian grid constitutes a limitation with respect to the spatial discretization, since it does not allow for adaptive refinement, for example in the vicinity of boundaries. This disadvantage can, however, be compensated by increasing the overall spatial resolution while the wall-clock time per solver step remains low.

Future work could investigate the weak pressure recovery observed downstream of the pressure minimum shown in the pressure coefficient analysis of the 3D sphere benchmark. Further studies could also consider substantially higher Reynolds numbers in order to test the algorithm's versatility and robustness. The analysis of additional flow quantities and more comparisons with high-resolution benchmarks might be necessary to expose and improve shortcomings of KSF-Chol and in order to make the algorithm fully suitable for scientific and engineering applications. The use of a more targeted sampling strategy for the SVD basis construction may improve approximation accuracy for a fixed number of sampled far-field offsets. A more sophisticated treatment of boundaries located in the computational domain could also be considered in order to improve the simulation accuracy of boundary-layer interactions. Finally, the computational methods developed here to accelerate kernel evaluations and kernel matrix inversions can be extended to a broader class of kernel-based approaches for solving PDEs, including those considered in \cite{ChenHosseiniOwhadiStuart2021, Bacho2026KROM, Bacho2023, BEM2023, BEM2019}.

\section*{Acknowledgements}
LS started the project under the Summer Undergraduate Research Fellowships (SURF) program at Caltech and continued working on it as part of an undergraduate thesis. LS thanks the SURF program and HO for funding, and AB and HO for their supervision and mentorship. AB gratefully acknowledges support from the German Research Foundation (DFG) through the Walter Benjamin Programme, as well as from the California Institute of Technology. 
HO acknowledges support from the Air Force Office of Scientific Research under MURI award number FOA-AFRL-AFOSR-2023-0004 (Mathematics of Digital Twins), the Department of Energy under award number DE-SC0023163 (SEA-CROGS: Scalable, Efficient, and Accelerated Causal Reasoning Operators, Graphs and Spikes for Earth and Embedded Systems) and the Vannevar Bush Faculty Fellowship under ONR-N000142512035.

\bibliographystyle{siamplain}
\bibliography{references}

@conference{StableFluid,
  author = {J. Stam},
  booktitle = {Proceedings of SIGGRAPH 1999},
  title = {Stable Fluids},
  year = {1999}
}

@article{PicTurb,
  author = {S. Elcott and Y. Tong and E. Kanso and P. Schr{\"o}der and M. Desbrun},
  title = {Stable, Circulation-Preserving, Simplicial Fluids},
  journal = {ACM Transactions on Graphics},
  year = {2007}
}

@article{sparse_chol_leibler,
  author = {F. Sch{\"a}fer and M. Katzfuss and H. Owhadi},
  title = {Sparse Cholesky Factorization by Kullback--Leibler Minimization},
  journal = {SIAM Journal on Scientific Computing},
  year = {2021}
}

@book{canuto1988spectral,
  author = {C. Canuto and M. Y. Hussaini and A. Quarteroni and T. A. Zang},
  title = {Spectral Methods in Fluid Dynamics},
  publisher = {Springer-Verlag},
  year = {1988}
}

@book{peyret2002spectral,
  author = {R. Peyret},
  title = {Spectral Methods for Incompressible Viscous Flow},
  publisher = {Springer},
  year = {2002}
}

@article{patterson1971spectral,
  author = {G. S. Patterson Jr. and S. A. Orszag},
  title = {Spectral Calculations of Isotropic Turbulence: Efficient Removal of Aliasing Interactions},
  journal = {Physics of Fluids},
  year = {1971}
}

@article{owhadi2023gph,
  author = {H. Owhadi},
  title = {Gaussian Process Hydrodynamics},
  journal = {Applied Mathematics and Mechanics (English Edition)},
  year = {2023}
}

@article{OberkampfTrucano2002,
  author = {W. L. Oberkampf and T. G. Trucano},
  title = {Verification and Validation in Computational Fluid Dynamics},
  journal = {Progress in Aerospace Sciences},
  year = {2002}
}

@article{msichol,
  author = {Y. Chen and H. Owhadi and F. Sch{\"a}fer},
  title = {Sparse Cholesky Factorization for Solving Nonlinear PDEs via Gaussian Processes},
  journal = {Mathematics of Computation},
  year = {2025}
}

@book{ferziger2020computational,
  author = {J. H. Ferziger and M. Peri{\'c} and R. L. Street},
  title = {Computational Methods for Fluid Dynamics},
  publisher = {Springer},
  edition = {4},
  year = {2020}
}

@article{orszag1971accuracy,
  author = {S. A. Orszag},
  title = {Numerical Simulation of Incompressible Flows within Simple Boundaries: Accuracy},
  journal = {Journal of Fluid Mechanics},
  year = {1971}
}

@article{MadayQuarteroni1982,
  author = {Y. Maday and A. Quarteroni},
  title = {Spectral and Pseudo-Spectral Approximations of the Navier--Stokes Equations},
  journal = {SIAM Journal on Numerical Analysis},
  year = {1982}
}

@article{Kraichnan1967,
  author = {R. H. Kraichnan},
  title = {Inertial Ranges in Two-Dimensional Turbulence},
  journal = {Physics of Fluids},
  year = {1967}
}

@book{Pope2000,
  author = {S. B. Pope},
  title = {Turbulent Flows},
  publisher = {Cambridge University Press},
  year = {2000}
}

@conference{MarefatAlamPope2022,
  author = {H. A. Marefat and J. Alam and K. Pope},
  booktitle = {Proceedings of the ASME 2022 Fluids Engineering Division Summer Meeting},
  title = {Immersed Boundary Method Implemented in LES for Flow Past a Sphere at Subcritical Reynolds Numbers},
  year = {2022}
}

@conference{ScholkopfHerbrichSmola2001,
  author = {B. Sch{\"o}lkopf and R. Herbrich and A. J. Smola},
  booktitle = {Computational Learning Theory},
  title = {A Generalized Representer Theorem},
  year = {2001}
}

@article{RakhubaOseledets2015,
  author = {M. V. Rakhuba and I. V. Oseledets},
  title = {Fast Multidimensional Convolution in Low-Rank Tensor Formats via Cross Approximation},
  journal = {SIAM Journal on Scientific Computing},
  year = {2015}
}

@book{ChanJin2007,
  author = {R. H.-F. Chan and X.-Q. Jin},
  title = {An Introduction to Iterative Toeplitz Solvers},
  publisher = {Society for Industrial and Applied Mathematics},
  year = {2007}
}

@article{Chorin1968,
  author = {A. J. Chorin},
  title = {Numerical Solution of the Navier--Stokes Equations},
  journal = {Mathematics of Computation},
  year = {1968}
}

@article{GuermondMinevShen2006,
  author = {J.-L. Guermond and P. Minev and J. Shen},
  title = {An Overview of Projection Methods for Incompressible Flows},
  journal = {Computer Methods in Applied Mechanics and Engineering},
  year = {2006}
}

@article{Costa2018,
  author = {P. Costa},
  title = {A FFT-Based Finite-Difference Solver for Massively-Parallel Direct Numerical Simulations of Turbulent Flows},
  journal = {Computers \& Mathematics with Applications},
  year = {2018}
}

@article{VanDerPoel2015,
  author = {E. P. van der Poel and R. Ostilla-M{\'o}nico and J. Donners and R. Verzicco},
  title = {A Pencil Distributed Finite Difference Code for Strongly Turbulent Wall-Bounded Flows},
  journal = {Computers \& Fluids},
  year = {2015}
}

@article{GreengardRokhlin1987,
  author = {L. Greengard and V. Rokhlin},
  title = {A Fast Algorithm for Particle Simulations},
  journal = {Journal of Computational Physics},
  volume = {73},
  number = {2},
  pages = {325--348},
  year = {1987},
  doi = {10.1016/0021-9991(87)90140-9}
}

@article{FongDarve2009,
  author = {W. Fong and E. Darve},
  title = {The Black-Box Fast Multipole Method},
  journal = {Journal of Computational Physics},
  volume = {228},
  number = {23},
  pages = {8712--8725},
  year = {2009},
  doi = {10.1016/j.jcp.2009.08.031}
}

@article{MessnerBramasCoulaudDarve2012,
  author = {M. Messner and B. Bramas and O. Coulaud and E. Darve},
  title = {Optimized M2L Kernels for the Chebyshev Interpolation Based Fast Multipole Method},
  journal = {arXiv preprint arXiv:1210.7292},
  year = {2012}
}

@book{Trefethen2013,
  author = {L. N. Trefethen},
  title = {Approximation Theory and Approximation Practice},
  publisher = {Society for Industrial and Applied Mathematics},
  year = {2013}
}

@article{EckartYoung1936,
  author = {C. Eckart and G. Young},
  title = {The Approximation of One Matrix by Another of Lower Rank},
  journal = {Psychometrika},
  volume = {1},
  number = {3},
  pages = {211--218},
  year = {1936},
  doi = {10.1007/BF02288367}
}

@article{BusingerGolub1965,
  author = {P. A. Businger and G. H. Golub},
  title = {Linear Least Squares Solutions by Householder Transformations},
  journal = {Numerische Mathematik},
  volume = {7},
  number = {3},
  pages = {269--276},
  year = {1965},
  doi = {10.1007/BF01436084}
}

@article{RodriguezEtAl2019,
  author = {I. Rodriguez and O. Lehmkuhl and M. Soria and S. G{\'o}mez and M. Dom{\'i}nguez-Pumar and L. Kowalski},
  title = {Fluid Dynamics and Heat Transfer in the Wake of a Sphere},
  journal = {International Journal of Heat and Fluid Flow},
  volume = {76},
  pages = {141--153},
  year = {2019},
  doi = {10.1016/j.ijheatfluidflow.2019.02.004}
}

@article{TomboulidesOrszag2000,
  author = {A. G. Tomboulides and S. A. Orszag},
  title = {Numerical Investigation of Transitional and Weak Turbulent Flow Past a Sphere},
  journal = {Journal of Fluid Mechanics},
  volume = {416},
  pages = {45--73},
  year = {2000},
  doi = {10.1017/S0022112000008880}
}

@article{NicoudDucros1999,
  author = {F. Nicoud and F. Ducros},
  title = {Subgrid-Scale Stress Modelling Based on the Square of the Velocity Gradient Tensor},
  journal = {Flow, Turbulence and Combustion},
  volume = {62},
  pages = {183--200},
  year = {1999},
  doi = {10.1023/A:1009995426001}
}

@article{BEM2023,
  author  = {Bacho, Aras},
  title   = {A Generalization of the {Moreau--Yosida} Regularization},
  journal = {Journal of Mathematical Analysis and Applications},
  volume  = {524},
  number  = {2},
  pages   = {127139},
  year    = {2023},
  doi     = {10.1016/j.jmaa.2023.127139},
  url     = {https://doi.org/10.1016/j.jmaa.2023.127139}
}

@article{BEM2019,
  author  = {Bacho, Aras and Emmrich, Etienne and Mielke, Alexander},
  title   = {An existence result and evolutionary
             {$\Gamma$}-convergence for perturbed gradient systems},
  journal = {Journal of Evolution Equations},
  volume  = {19},
  pages   = {479--522},
  year    = {2019},
  doi     = {10.1007/s00028-019-00484-x},
  url     = {https://doi.org/10.1007/s00028-019-00484-x}
}

@article{Bacho2026KROM, 
  author        = {Bacho, Aras and Lee, Jonghyeon
                   and Owhadi, Houman},
  title         = {{KROM}: Kernelized Reduced Order Modeling},
  year          = {2026},
  eprint        = {2603.00360},
  archivePrefix = {arXiv},
  primaryClass  = {math.NA},
  doi           = {10.48550/arXiv.2603.00360},
  url           = {https://arxiv.org/abs/2603.00360}
}

@article{Bacho2023,
  author  = {Bacho, Aras},
  title   = {Abstract nonlinear evolution inclusions of
             second order with applications in
             visco-elasto-plasticity},
  journal = {Journal of Differential Equations},
  volume  = {363},
  pages   = {126--169},
  year    = {2023},
  doi     = {10.1016/j.jde.2023.03.008},
  url     = {https://doi.org/10.1016/j.jde.2023.03.008}
}

@article{Bacho2026Operator,
  author  = {Bacho, Aras and Sorokin, Aleksei G.
             and Yang, Xianjin and Bourdais, Th{\'e}o
             and Calvello, Edoardo and Darcy, Matthieu
             and Hsu, Alexander W. and Hosseini, Bamdad
             and Owhadi, Houman},
  title   = {Operator Learning at Machine Precision},
  journal = {Journal of Computational Physics},
  volume  = {566},
  pages   = {115240},
  year    = {2026},
  doi     = {10.1016/j.jcp.2026.115240},
  url     = {https://doi.org/10.1016/j.jcp.2026.115240}
}

@article{OwhadiScovelSchaefer2019,
  author = {H. Owhadi and C. Scovel and F. Sch{\"a}fer},
  title = {Statistical Numerical Approximation},
  journal = {Notices of the American Mathematical Society},
  volume = {66},
  number = {10},
  pages = {1608--1617},
  year = {2019},
  doi = {10.1090/noti1963}
}

@article{ChenHosseiniOwhadiStuart2021,
  author = {Y. Chen and B. Hosseini and H. Owhadi and A. M. Stuart},
  title = {Solving and Learning Nonlinear {PDEs} with Gaussian Processes},
  journal = {Journal of Computational Physics},
  volume = {447},
  pages = {110668},
  year = {2021},
  doi = {10.1016/j.jcp.2021.110668}
}

\appendix

\section{Numerical parameters and stability criteria}

\subsubsection{Kernel choice, preconditioning, and sparsity parameter}
Two different kernel families are considered in the numerical experiments.
For the energy-cascade and merging-vortices test cases, we use a Gaussian
kernel whose length scale $\ell$ scales with the spatial grid spacing $dx=\frac{L_{x}}{N_{x}}$ and $dy$. For the von K\'arm\'an benchmark, we use a Mat\'ern kernel with smoothness $\nu_k=5/2$, length scale $\ell=0.22$, and amplitude $\sigma_K=0.4$. For the sphere benchmark, we use $\nu_K=5/2$, $\ell=0.12$, and $\sigma_K=1.0$. The parameters were selected empirically. The Gaussian kernel produced non-physical oscillatory artifacts in the higher-Reynolds-number test cases.

A regularization nugget $\varepsilon_{\mathrm{nug}}=10^{-8}$ is added to the
diagonal of the kernel matrix in order to improve its conditioning and to stabilize the subsequent linear solves. 

The sparsity of the Cholesky factorization is controlled by
$\rho_{\mathrm{chol}}$. All benchmarks use $\rho_{\mathrm{chol}}=7.5$, independently of grid size and spatial dimension. In contrast, the approximation accuracy result of \cite{sparse_chol_leibler} suggests $\rho_{\mathrm{chol}}\approx\log(N/\varepsilon_{\mathrm{chol}})$. With a fixed sparsity parameter, we observe increasing PCG iteration numbers on larger grids, although the numbers remain constant during each simulation. We attribute this to a less accurate inverse kernel approximation $\widehat{\mathbf{K}}^{-1}\approx\mathbf{L}\mathbf{L}^T$ for a fixed $\rho_{\mathrm{chol}}$.

For two-dimensional problems, \cite{msichol} suggests a substantially smaller value of $\rho_{\mathrm{chol}}=3.2$. In our higher-Reynolds-number simulations, however, this choice did not lead to sufficiently small iteration numbers in the preconditioned iterative solves. With $\rho_{\mathrm{chol}}=7.5$, the number of iterations typically remains in the single-digit range and can increase as the flow becomes more turbulent.

\section{Benchmark Setup Details and Results}
\label{app:benchmark-setups}

\subsection{Energy Spectrum Benchmark}
\label{app:energy-spectrum-setup}

\paragraph{Forcing}
In this benchmark, we include an external forcing term consisting of two components. The first is a forcing acting on a rectangular subdomain:
\begin{equation}
f_{\text{local}}(x,y) = \frac{2\pi}{\ell_x} \cos\!\left(\frac{2\pi(x-0.4L_x)}{\ell_x}\right),
\quad \ell_x = 0.2L_x,
\end{equation}
restricted to $x \in [0.4L_x,0.6L_x]$ and $y \in [0.1L_y,0.3L_y]$.
In addition, we apply a smooth global forcing,
\begin{equation}
f_{\text{global}}(x,y) = 0.05\bigl(\sin(k_f x + k_f y) + \cos(k_f x + k_f y)\bigr), 
\quad k_f = 1.
\end{equation}
The global forcing excites large-scale Fourier modes associated with the domain size, which is essential to sustain energy at low wavenumbers and obtain a physically meaningful energy spectrum.

\paragraph{Stable Fluids}
We adapt the Stable Fluids algorithm to the vorticity formulation of the NS-equation.

\paragraph{KSF-Chol and Kernel choice}
\phantomsection
\label{par:energy-ksf-kernel}
We implement a periodic Gaussian kernel obtained via embedding the torus into a higher-dimensional Euclidean space. Specifically, each point $\mathbf{x} = (x,y)$ is mapped to
\begin{equation}
\Phi(\mathbf{x}) = \bigl(\cos(2\pi x/L_x), \sin(2\pi x/L_x), \cos(2\pi y/L_y), \sin(2\pi y/L_y)\bigr),
\end{equation}
which enforces periodicity by construction. The kernel is defined as
\begin{equation}
G(\mathbf{x},\mathbf{y}) 
= \exp\!\left(-\frac{1}{2\ell^2} \|\Phi(\mathbf{x}) - \Phi(\mathbf{y})\|^2 \right).
\end{equation}
More parameters of the simulation are summarized in \autoref{tab:vortex_pair_params}.

\begin{table}[h]
\centering
\caption{Parameters for the energy spectrum setup.}
\label{tab:vortex_pair_params}
\scriptsize
\begin{tabular}{@{}ll@{}}
\toprule
\textbf{Solver} & \textbf{KSF-Chol} \\ \midrule
Domain size           & $2\pi \times 2\pi$                                   \\
$N_x, N_y$            & $150, 150$                                           \\
Viscosity $\nu$       & $0.001$                                              \\
Forcing amplitude     & $2.0$                                                \\
Timestep              & $0.01$                                               \\
Kernel type           & Gaussian                                             \\
Kernel length scale $\ell$ & $0.04$                                          \\
sparsity $\rho_{\mathrm{chol}}$      & $7.5$                                                \\
 \bottomrule
\end{tabular}
\end{table}

\subsection{Two Merging Vortices Benchmark}
\label{app:merging-vortices-setup}

\paragraph{Setup}
The setup is largely identical to the one described in the previous section concerning the energy cascade experiment, except that the computational domain is now $\mathbb{T}^2 = [0,\pi]^2$ and no external forcing is applied. We again consider the vorticity--streamfunction formulation of the incompressible NS equations, and the KSF-Chol solver remains unchanged. The initial condition consists of two Gaussian vortices and is shown in \autoref{fig:merging_vort_pic} (left, $t=0$). We use the same periodic Gaussian kernel constructed via embedding.

\begin{table}[h]
\centering
\caption{Parameters for the merging vortices setup.}
\label{tab:merging_vort_tab}
\scriptsize
\setlength{\tabcolsep}{3pt}
\begin{tabularx}{\linewidth}{@{}>{\raggedright\arraybackslash}X *{3}{>{\centering\arraybackslash}X}@{}}
\toprule
\textbf{Solver} & \textbf{Pseudospectral} & \textbf{Stable Fluids} & \textbf{KSF-Chol} \\ \midrule
Domain size& $\pi \times \pi$ & $\pi \times \pi$ & $\pi \times \pi$ \\
Resolution & $120^2$ & $120^2$ & $120^2$ \\
Viscosity & $10^{-3}$ & $10^{-3}$ & $10^{-3}$ \\
Timestep & $0.001$ & $0.01$ & $0.01$ \\
Wall-clock time per time step & $0.010171 \mathrm{s}$ & $0.014410 \mathrm{s}$ & $0.272236 \mathrm{s}$ \\
Kernel & --- & --- & Gaussian \\
$\ell$ & --- & --- & $0.05$ \\
$\rho_{\mathrm{chol}}$ & --- & --- & $7.5$ \\ \bottomrule
\end{tabularx}
\end{table}

\subsection{Von Kármán Vortex Street Benchmark}
\label{app:karman-setup}

\paragraph{Boundary conditions}
No external forcing is applied. The velocity boundary conditions are prescribed as
\begin{equation}
\begin{aligned}
u &= U_{\mathrm{in}}, \quad v=0,
&& x=0,\\
\partial_x u &= 0, \quad \partial_x v=0,
&& x=L_x,\\
u &= 0, \quad v=0,
&& y\in\{0,L_y\},\\
u &= 0, \quad v=0,
&& \mathbf{x}\in\partial\Omega_{\mathrm{cyl}}.
\end{aligned}
\end{equation}

\paragraph{Stable Fluids}

We use the Stable Fluids method in the velocity
formulation with the external
forcing set to zero, $\mathbf{f}=0$.

\paragraph{KSF-Chol and Kernel choice}
\phantomsection
\label{par:karman-ksf-kernel}
For this benchmark we use the KSF-Chol method with a Matérn kernel and smoothness parameter $\nu_K = 5/2$:
\begin{equation}
K(r) = \sigma_K^2 \left( 1 + \frac{\sqrt{5}\,r}{\ell} + \frac{5 r^2}{3\ell^2} \right)
\exp\left(-\frac{\sqrt{5}\,r}{\ell}\right),
\end{equation}
where $r = \|\mathbf{x}-\mathbf{x}'\|$. A Gaussian kernel was also tested but led to numerical artifacts. 

\begin{table}[h]
\centering
\caption{Parameters for the von Kármán vortex street setup.}
\label{tab:karman_setup}
\scriptsize
\setlength{\tabcolsep}{3pt}
\begin{tabularx}{\linewidth}{@{}>{\raggedright\arraybackslash}X *{4}{>{\centering\arraybackslash}X}@{}}
\toprule
\textbf{Solver} & \textbf{Reference} & \textbf{Stable Fluids} & \textbf{KSF-Chol} & \textbf{KSF-Chol}\\
 & & & (dense) & (approx) \\
\midrule
Domain size & $4\times2$ & $4\times2$ & $4\times2$ & $4\times2$ \\
Resolution & $180\times90$ & $180\times90$ & $180\times90$ & $180\times90$ \\
Viscosity & $0.002$ & $0.002$ & $0.002$ & $0.002$ \\
Timestep & $0.001$ & $0.001$ & $0.001$ & $0.001$ \\
Kernel & --- & --- & Matérn & Matérn \\
$\sigma_K$ & --- & --- & $0.4$ & $0.4$ \\
$\ell$ & --- & --- & $0.22$ & $0.22$ \\
$\rho_{\mathrm{chol}}$ & --- & --- & $7.5$ & $7.5$ \\
\bottomrule
\end{tabularx}
\end{table}

\begin{figure}[h]
\centering
\includegraphics[width=0.75\linewidth]{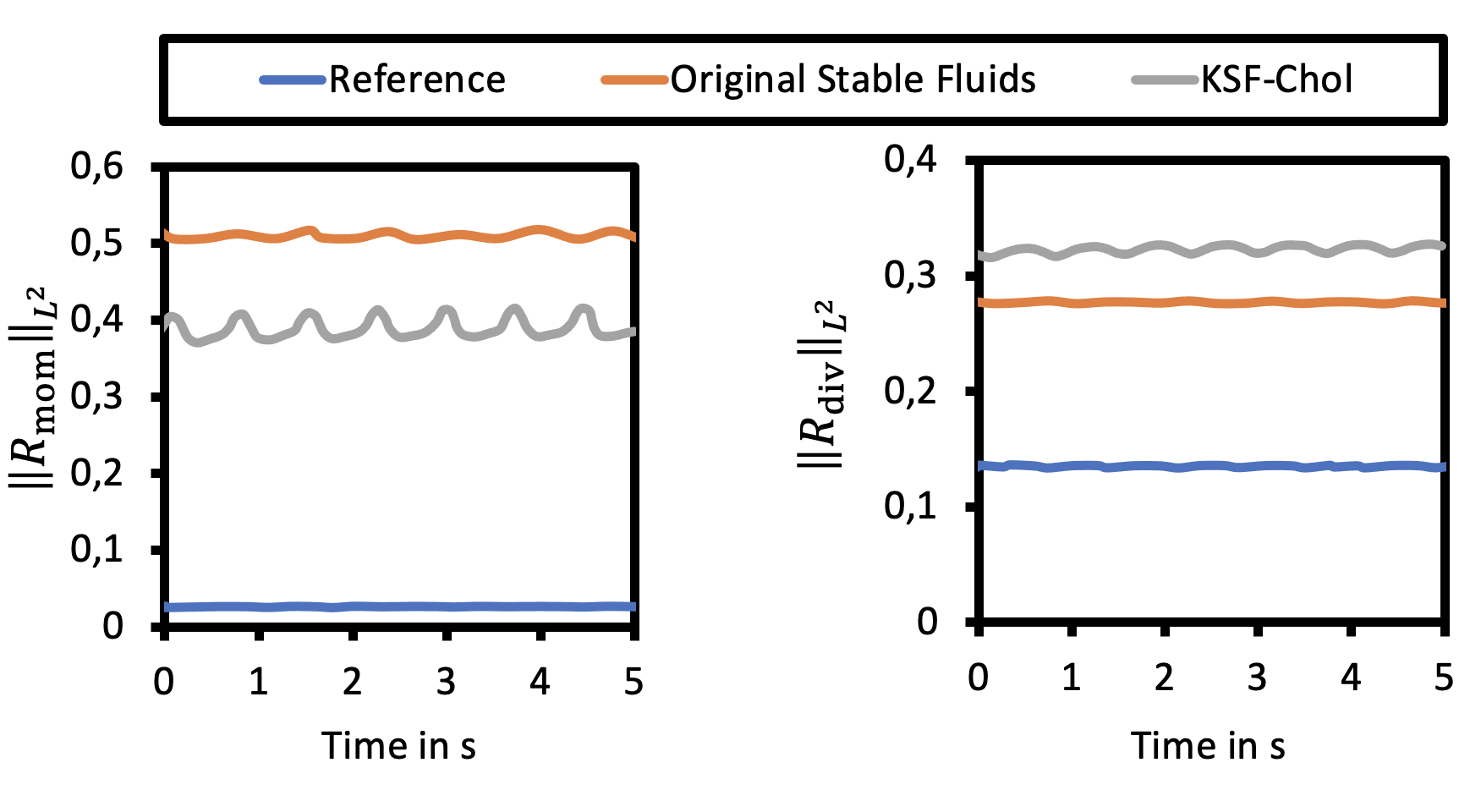}
\caption[]{\textbf{Left:} Momentum residual in $L^2$ norm of von Kármán vortex street simulation. \textbf{Right:} Divergence residual in $L^2$ norm of von Kármán vortex street simulation.}
\label{fig:mom_res_karman}
\end{figure}

\subsection{3D Flow Past a Sphere at $\mathrm{Re}=1000$}
\label{app:sphere-setup}

\paragraph{Boundary conditions}
The velocity boundary conditions are prescribed as
\begin{equation}
\begin{aligned}
\mathbf{u} &= (U_{\mathrm{in}},0,0)^{\mathsf T},
&& x=0,\\
\partial_x \mathbf{u} &= \mathbf{0},
&& x=L_x,\\
v=0,\quad \partial_y u=\partial_y w &= 0,
&& y\in\{0,L_y\},\\
w=0,\quad \partial_z u=\partial_z v &= 0,
&& z\in\{0,L_z\}.
\end{aligned}
\end{equation}

A stationary spherical obstacle with center $\mathbf{x}_{\mathrm{s}}$ and
radius $R_{\mathrm{s}}$ is applied according to
$
\Omega_{\mathrm{s}}
=
\bigl\{
\mathbf{x}\in\Omega :
\|\mathbf{x}-\mathbf{x}_{\mathrm{s}}\|_2
\le R_{\mathrm{s}}
\bigr\}.
$
The no-slip condition on the spherical obstacle is prescribed as
$
\mathbf{u}(\mathbf{x},t)=\mathbf{0}
$ with $
\mathbf{x}\in\Omega_{\mathrm{s}}.
$

\begin{table}[t]
\centering
\begin{minipage}[t]{0.48\linewidth}
\vspace{0pt}
\centering
\caption{Parameters for the 3D flow over a sphere at $Re=1000$.}
\label{tab:setup_re1000}
\scriptsize
\setlength{\tabcolsep}{3pt}
\begin{tabularx}{\linewidth}{@{}>{\raggedright\arraybackslash}X >{\raggedright\arraybackslash}p{0.32\linewidth}@{}}
\toprule
\textbf{Parameter} & \textbf{Value} \\
\midrule
Domain size                         & $4.0 \times 2.0 \times 2.0$ \\
$N_x, N_y, N_z$                     & $180,90,90$ \\
Sphere center $\mathbf{x}_{\mathrm{s}}$ & $(0.7,1.0,1.0)$ \\
Sphere radius $R_{\mathrm{s}}$        & $0.17$ \\
Viscosity $\nu$                     & $0.00034$ \\
$U_{\mathrm{inlet}}$                & $1.0$ \\
Forcing amplitude                   & $0.0$ \\
Timestep $\Delta t$                 & $0.01$ \\
Kernel type                         & Matérn-$5/2$ \\
Kernel amplitude $\sigma_K$           & $1.0$ \\
Kernel length scale $\ell$             & $0.12$ \\
sparsity $\rho_{\mathrm{chol}}$              & $7.5$ \\
\bottomrule
\end{tabularx}
\end{minipage}
\hfill
\begin{minipage}[t]{0.48\linewidth}
\vspace{0pt}
\centering
\caption{Parameters of the near-/far-field kernel evaluation for the 3D flow over a sphere benchmark.}
\label{tab:cheb_svd_fft_re1000}
\scriptsize
\setlength{\tabcolsep}{3pt}
\begin{tabularx}{\linewidth}{@{}>{\raggedright\arraybackslash}X >{\raggedright\arraybackslash}p{0.32\linewidth}@{}}
\toprule
\textbf{Parameter} & \textbf{Value} \\
\midrule
Near-field radius $n_{\mathrm{near}}$ & $3$ cells \\
Near-field stencil                    & $7^3$ cells \\
Chebyshev order $p_{\mathrm{cheb}}$                    & $10$ \\
Chebyshev--Lobatto nodes $L$           & $1331$ \\
Total-degree coefficients $P$          & $286$ \\
Total far-field offsets $|\mathcal{F}_c|$ & $11\,502\,376$ \\
Number of sampled far-field offsets $S$ & $8192$ \\
SVD tolerance                          & $10^{-10}$ \\
Maximum SVD rank $R$                       & $64$ \\
\midrule
Dense validation interval              & $100$ steps \\
Dense evaluation technique             & PyKeOps \\
\bottomrule
\end{tabularx}
\end{minipage}
\end{table}

\section{Reference Solvers}
\label{app:reference-solvers}

\subsection{Fourier Pseudo-Spectral Reference Solver}
\label{app:pseudospectral-reference-solver}
The pseudo-spectral solver as shown in Algorithm~\ref{alg:dns} advances the two-dimensional NS equations
in vorticity form on a periodic domain.
The Fourier wave numbers are denoted by $k_x$ and $k_y$, and
\begin{equation}
|\mathbf{k}|^2 = k_x^2 + k_y^2.
\end{equation}
Spatial derivatives are evaluated directly in Fourier space. In particular,
differentiation with respect to $x$ and $y$ corresponds to multiplication with
$i k_x$ and $i k_y$. The Laplacian is diagonal in Fourier
space and is represented by multiplication with $-|\mathbf{k}|^2$
\cite{peyret2002spectral,canuto1988spectral}. To suppress aliasing, a spectral
mask is introduced,
\begin{equation}
\chi(\mathbf{k}) =
\begin{cases}
1,
&
|k_x| < \dfrac{2}{3}k_{x,\max}
\ \text{and}
|k_y| < \dfrac{2}{3}k_{y,\max},
\\[0.5em]
0,
&
\text{otherwise}.
\end{cases}
\end{equation}
High-frequency Fourier modes outside the
$2/3$ range are therefore removed
\cite{patterson1971spectral}.
For a given spectral vorticity field, the velocity is reconstructed through the
streamfunction $\psi$. The relation
$
-\Delta \psi = \omega
$
becomes
$
|\mathbf{k}|^2\hat{\psi}
=
\hat{\omega}
$
in Fourier space. Hence,
$
\hat{\psi}
=
\frac{\hat{\omega}}{|\mathbf{k}|^2},
$ and $
\hat{\psi}(\mathbf{0})=0,
$
and the velocity components are obtained as
$
u
=
\mathcal{F}^{-1}
\left(
i k_y \hat{\psi}
\right),
$ and $
v
=
\mathcal{F}^{-1}
\left(
-i k_x \hat{\psi}
\right).
$
The spatial derivatives of the vorticity are also evaluated spectrally by 
$
\omega_x
=
\mathcal{F}^{-1}
\left(
i k_x\hat{\omega}
\right),
$ and $
\omega_y
=
\mathcal{F}^{-1}
\left(
i k_y\hat{\omega}
\right).
$
The nonlinear advection term is formed in physical space as
$
\mathcal{N}(\omega)
=
u\,\omega_x
+
v\,\omega_y
$
and is subsequently transformed back to spectral space using the dealiased
transform 
$
\widehat{\mathcal{N}}
=
\mathcal{F}_d
\left(
\mathcal{N}(\omega)
\right).
$
The resulting spectral vorticity equation is
\begin{equation}
\frac{\partial\hat{\omega}}{\partial t}
=
-\widehat{\mathcal{N}}
-
\nu |\mathbf{k}|^2\hat{\omega}.
\end{equation}
For compactness, its right-hand side is denoted by
$
\mathcal{R}(\hat{\omega})
$

Limitations are the restriction to simple
periodic geometries, the reduced effectiveness of Fourier approximations
for nonsmooth solutions, and the time-step restriction associated with the stability of the
explicit RK4 treatment of both advection and diffusion.

\begin{algorithm}
\caption{Pseudo-Spectral Navier--Stokes Solver in Vorticity Formulation}
\label{alg:dns}

\normalsize
\algrenewcommand\alglinenumber[1]{\normalsize #1}

\begin{algorithmic}[1]

\State \textbf{Input:} spatial resolution $N_x \times N_y$, domain
$[0,L_x]\times[0,L_y]$, time step $\Delta t$, number of time steps
$N_t$, viscosity $\nu$, and initial vorticity field $\omega^0(\mathbf{x})$

\vspace{0.3em}
\Statex \textbf{Construct the $2/3$ dealiasing filter:}

\State Define the dealiased Fourier transform
$\mathcal{F}_d(f)
=
\chi(\mathbf{k})\,\mathcal{F}(f)$

\vspace{0.3em}

\State Transform the initial vorticity:
$\hat{\omega}^{\,0}
=
\mathcal{F}_d
\left(
\omega^0
\right)$

\vspace{0.3em}

\For{$n=0,1,\dots,N_t-1$}

\Statex \textbf{Advance the spectral vorticity using RK4:}

\State
$\mathbf{k}_1
=
\mathcal{R}
\left(
\hat{\omega}^{\,n}
\right)$

\State
$\mathbf{k}_2
=
\mathcal{R}
\left(
\hat{\omega}^{\,n}
+
\frac{\Delta t}{2}\mathbf{k}_1
\right)$

\State
$\mathbf{k}_3
=
\mathcal{R}
\left(
\hat{\omega}^{\,n}
+
\frac{\Delta t}{2}\mathbf{k}_2
\right)$

\State
$\mathbf{k}_4
=
\mathcal{R}
\left(
\hat{\omega}^{\,n}
+
\Delta t\,\mathbf{k}_3
\right)$

\State Update the spectral vorticity:
$\hat{\omega}^{\,n+1}
=
\hat{\omega}^{\,n}
+
\frac{\Delta t}{6}
\left(
\mathbf{k}_1
+
2\mathbf{k}_2
+
2\mathbf{k}_3
+
\mathbf{k}_4
\right)$

\State Recover the physical-space vorticity:
$\omega^{n+1}
=
\mathcal{F}^{-1}
\left(
\hat{\omega}^{\,n+1}
\right)$

\EndFor

\vspace{0.3em}

\State \textbf{Return:} sequence
$\{\omega^n\}_{n=0}^{N_t}$

\end{algorithmic}
\end{algorithm}

\subsection{Implementation Details of the Finite-Difference Reference Solver}
\label{app:fd-reference-solver}

For the pressure projection, the discrete divergence of the intermediate
velocity is evaluated at interior grid points by
\begin{equation}
\left(
\nabla_h\cdot\mathbf{u}^{*}
\right)_{i,j}
=
\frac{
u^{*}_{i,j+1}
-
u^{*}_{i,j-1}
}{
2\Delta x
}
+
\frac{
v^{*}_{i+1,j}
-
v^{*}_{i-1,j}
}{
2\Delta y
}.
\end{equation}
The resulting pressure-Poisson problem is diagonalized by a
two-dimensional discrete cosine transform (DCT).
The complete algorithm is summarized in ~\ref{alg:fd_projection}.

\begin{algorithm}
\caption{Finite-Difference Projection Solver for Flow Past a Cylinder}
\label{alg:fd_projection}

\normalsize
\algrenewcommand\alglinenumber[1]{\normalsize #1}

\begin{algorithmic}[1]

\State \textbf{Input:} spatial resolution $N_x\times N_y$, domain size $\Omega$,
time step $\Delta t$, number of time steps $N_t$, inlet velocity
$U_{\mathrm{in}}$, density $\rho_f$, Reynolds number $\mathrm{Re}$, and cylinder
geometry

\vspace{0.3em}

\State Initialize $u^0$, $v^0$, and $p^0$

\vspace{0.3em}

\For{$n=0,1,\dots,N_t-1$}

    \Statex \textbf{Compute intermediate velocity:}

    \State $\displaystyle
    u^{\ast}
    =
    u^n
    -
    \Delta t
    \left(
    u^n D_xu^n
    +
    v^n D_yu^n
    \right)
    +
    \nu\Delta t\,\Delta_hu^n$

    \State $\displaystyle
    v^{\ast}
    =
    v^n
    -
    \Delta t
    \left(
    u^n D_xv^n
    +
    v^n D_yv^n
    \right)
    +
    \nu\Delta t\,\Delta_hv^n$

    \Statex \textbf{Enforce boundary condition inside the cylinder:}

    \State $\mathbf{u}^{\ast}(\mathbf{x})=\mathbf{0},
    \qquad
    \mathbf{x}\in\Omega_{\mathrm{cyl}}$

    \vspace{0.3em}

    \Statex \textbf{Solve pressure Poisson equation:}

    \State $\displaystyle
    \Delta p^{n+1}
    =
    \frac{\rho_f}{\Delta t}
    \left(
    D_xu^{\ast}
    +
    D_yv^{\ast}
    \right)$

    \State Transform Poisson equation using the two-dimensional DCT

    \State $\displaystyle
    \widehat{p}_{k,l}
    =
    \frac{\widehat{b}_{k,l}}{\lambda_{k,l}},
    \qquad
    \widehat{p}_{0,0}=0$

    \State Recover $p^{n+1}$ using the inverse DCT

    \vspace{0.3em}

    \Statex \textbf{Project velocity field:}

    \State $\displaystyle
    u^{n+1}
    =
    u^{\ast}
    -
    \frac{\Delta t}{\rho_f}D_xp^{n+1}$

    \State $\displaystyle
    v^{n+1}
    =
    v^{\ast}
    -
    \frac{\Delta t}{\rho_f}D_yp^{n+1}$

    \Statex \textbf{Apply boundary conditions:}

    \State $\displaystyle
    u^{n+1}=U_{\mathrm{in}},
    \quad
    v^{n+1}=0
    \qquad
    \text{on } \Gamma_{\mathrm{in}}$

    \State $\displaystyle
    u^{n+1}=0,
    \quad
    v^{n+1}=0
    \qquad
    \text{on } \Gamma_{\mathrm{wall}}\cup\Omega_{\mathrm{cyl}}$

    \State $\displaystyle
    \frac{\partial u^{n+1}}{\partial x}
    =
    \frac{\partial v^{n+1}}{\partial x}
    =
    \frac{\partial p^{n+1}}{\partial x}
    =0
    \qquad
    \text{on } \Gamma_{\mathrm{out}}$

\EndFor

\vspace{0.3em}

\State \textbf{Return:} vorticity fields
$\{\omega^{n_s}\}_{s=0}^{N_s}$

\end{algorithmic}
\end{algorithm}

\newpage


\end{document}